\documentclass[10pt,a4paper]{article}
\usepackage{bbm}
\usepackage{mathrsfs}
\usepackage{amsfonts}
\usepackage{}

\usepackage{tikz}
\usetikzlibrary{positioning}
\usetikzlibrary{arrows}

\usepackage[leqno]{amsmath}
\usepackage{graphicx}
\usepackage{latexsym}
\usepackage{amsmath,amsfonts,amssymb,amsthm,mathrsfs,euscript,makeidx,color}
\usepackage{enumerate}
\usepackage{color}
\usepackage[colorlinks,linkcolor=blue,anchorcolor=green,citecolor=red]{hyperref}
\usepackage{indentfirst}
\renewcommand{\baselinestretch}{1.0}
\def\sqr#1#2{{\vcenter{\vbox{\hrule height.#2pt
              \hbox{\vrule width.#2pt height#1pt \kern#1pt \vrule width.#2pt}
              \hrule height.#2pt}}}}
\def\signed #1{{\unskip\nobreak\hfil\penalty50
              \hskip2em\hbox{}\nobreak\hfil#1
              \parfillskip=0pt \finalhyphendemerits=0 \par}}
\def\endpf{\signed {$\sqr69$}}

\def\5n{\negthinspace \negthinspace \negthinspace \negthinspace \negthinspace }
\def\4n{\negthinspace \negthinspace \negthinspace \negthinspace }
\def\3n{\negthinspace \negthinspace \negthinspace }
\def\2n{\negthinspace \negthinspace }
\def\1n{\negthinspace }

\def\dbH{\mathbb{H}}

\def\dbR{\mathbb{R}}

\def\sD{\mathscr{D}}

\def\sK{\mathscr{K}}

\def\sO{\mathscr{O}}

\def\sY{\mathscr{Y}}

\def\cK{{\cal K}}

\def\cl{{\cal l}}

\def\BN{{\bf N}}

\def\co{{\mathop{\rm co}\,}}
\def\ds{\displaystyle}

\def\ns{\noalign{\ss}}
\def\no{\noindent}

\def\ss{\smallskip}
\def\ms{\medskip}

\def\q{\quad}
\def\qq{\qquad}
\def\hb{\hbox}

\def\({\Big (}
\def\){\Big )}
\def\[{\Big[}
\def\]{\Big]}

\def\lan{\langle}
\def\ran{\rangle}

\def\rf{\eqref}

\def\a{\alpha}
\def\b{\beta}

\def\d{\delta}
\def\e{\varepsilon}
\def\z{\zeta}

\def\l{\lambda}
\def\m{\mu}
\def\n{\nu}

\def\t{\tau}
\def\f{\varphi}

\def\i{\infty}

\def\Ra{\mathop{\Rightarrow}}
\def\La{\mathop{\Leftarrow}}

\def\h{\widehat}
\def\wt{\widetilde}

\def\cd{\cdot}
\def\cds{\cdots}

\def\cl{\overline}
\def\pa{\partial}
\def\deq{\triangleq}
\def\les{\leqslant}
\def\ges{\geqslant}
\def\bde{\begin{definition}\label}
\def\ede{\end{definition}}
\def\be{\begin{equation}}
\def\bel{\begin{equation}\label}
\def\ee{\end{equation}}
\def\bt{\begin{theorem}\label}
\def\et{\end{theorem}}
\def\bc{\begin{corollary}\label}
\def\ec{\end{corollary}}
\def\bl{\begin{lemma}\label}
\def\el{\end{lemma}}
\def\bp{\begin{proposition}\label}
\def\ep{\end{proposition}}
\def\bas{\begin{assumption}\label}
\def\eas{\end{assumption}}
\def\br{\begin{remark}\label}
\def\er{\end{remark}}
\def\bex{\begin{example}\label}
\def\ex{\end{example}}
\def\ba{\begin{array}}
\def\ea{\end{array}}
\def\ben{\begin{enumerate}}
\def\een{\end{enumerate}}

\def\square#1{\vbox{\hrule\hbox{\vrule height#1%
     \kern#1\vrule}\hrule}}
\def\rectangle#1#2{\vbox{\hrule\hbox{\vrule height#1%
     \kern#2\vrule}\hrule}}

\font\tenbb=msbm10 \font\sevenbb=msbm7 \font\fivebb=msbm5

\newfam\bbfam
\scriptscriptfont\bbfam=\fivebb \textfont\bbfam=\tenbb
\scriptfont\bbfam=\sevenbb

\newtheorem{theorem}{\indent Theorem}[section]
\newtheorem{definition}[theorem]{\indent Definition}
\newtheorem{proposition}[theorem]{\indent Proposition}
\newtheorem{corollary}[theorem]{\indent Corollary}
\newtheorem{lemma}[theorem]{\indent Lemma}
\newtheorem{remark}[theorem]{\indent Remark}
\newtheorem{example}[theorem]{\indent Example}

\newtheorem{assumption}[theorem]{\indent Assumption}

\makeatletter
   
   \@addtoreset{equation}{section}
\makeatother

\begin{document}

\title{\bf Continuity of the Value Function for Deterministic Optimal Impulse Control with Terminal State Constraint}

\author{Yue Zhou\footnote{School of Mathematics and Statistics, Central South University, Changsha, Hunan, China. Email:{\tt zhouyuemath@csu.edu.cn.}},
~~~Xinwei Feng\footnote{Zhongtai Securities Institute for Financial Studies, Shandong University, Jinan, Shandong
250100, China. Email:{\tt xwfeng} {\tt @sdu.edu.cn.}},~~~Jiongmin Yong\footnote{Department of
Mathematics, University of Central Florida, Orlando, FL 32816, USA. This author was supported in part by NSF Grant DMS-1812921. Email:{\tt jiongmin.yong@ucf.edu.} }}

\maketitle

\no\bf Abstract: \rm Deterministic optimal impulse control problem with terminal state constraint is considered. Due to the appearance of the terminal state constraint, the value function might be discontinuous in general. The main contribution of this paper is the introduction of an intrinsic condition under which the value function is continuous. Then by a Bellman dynamic programming method, the corresponding Hamilton-Jacobi-Bellman type quasi-variational inequality (QVI, for short) is derived for which the value function is a viscosity solution. The issue of whether the value function is characterized as the unique viscosity solution to this QVI is carefully addressed and the answer is left open challengingly.

\ms

\no\bf Keywords: \rm Optimal impulse control, terminal state constraint, continuity of value function, Hamilton-Jacobi-Bellman quasi-variational inequality, viscosity solution.

\ms

\no\bf AMS Mathematics Subject Classification. \rm 49N25, 49L20, 49L25.

\section{Introduction}

It is well-known that in general classical continuous-time optimal control theory, there are two major approaches: variational method leading to Pontryagin's maximum principle (MP, for short), and dynamic programming method leading to Hamilton-Jacobi-Bellman (HJB, for short) equation. The former works for the problems with possible terminal state constraint and it gives necessary conditions for (possibly existed) open-loop optimal controls (\cite{Pontryagin 1962, Yong-Zhou 1999}). The latter works for the problems {\sl without} terminal state constraint and it leads to a characterization of the value function as the unique viscosity solution to the HJB equation, then formally, optimal control of state feedback form can be obtained (\cite{Bellman 1957, Fleming-Soner 1992, Li-Yong 1995, Bardi-Capuzzo-Dolcetta 1997, Yong-Zhou 1999}). In the case that there exists a terminal state constraint, the value function might not be everywhere defined (which involves the controllability issue), and even in the set on which the value function is defined, due to the set of admissible controls is depending on the initial pair, the continuity of the value function is not guaranteed. As a result, the corresponding satisfactory general theory of viscosity solution to the HJB equation for the terminal state constrained is not available as of today. Therefore, people had made some efforts to introduce proper conditions so that the continuity of the value function can still be obtained. One such an effort is for time optimal control problem (with a target set which is a terminal state constraint). To ensure the continuity of the value function, the so-called {\it small time local controllability} (STLC, for short) was introduced by Sussmann in 1987 (\cite{Sussmann 1987, Bardi-Capuzzo-Dolcetta 1997, Yong 2015}). This condition means that when the state gets close to the boundary of the target set (from outside), only a small amount of time is needed to drive the state to the target by a control action. This then will lead to the continuity of the value function.

\ms

For optimal impulse control problems, similar to the continuous control case, one also has two major approaches. For variational method, there are works on MP; we mention  \cite{Rempala-Zabczyk 1988,  Hu-Yong 1991, Yong-Zhang 1992, Dykhta-Samsonuk 2001, Chahim-Hartl-Kort 2012, Wu-Zhang 2012, Feng 2016}, for a partial list. On the other hand, since the initiation of optimal impulse control problems by Bensoussan--Lions in the early 1970s (\cite{Bensoussan-Lions 1973, Bensoussan-Lions 1984}), the dynamic programming method has been a very popular approach to the problem. It is interesting that the corresponding HJB equation is a {\it quasi-variational inequality} (QVI, for short) to which the value function is the unique viscosity solution (\cite{Barles 1985, Barles 1985b}) under proper conditions. There are quite a few follow-up works, see \cite{Menaldi 1987, Lenhart 1989, Yong 1989, Tang-Yong 1993, Yong 1994, Li-Yong 1995, Menaldi-Robin 2017, Korn-Melnyk-Seifried 2017, Belak-Christensen-Seifried 2019} for a partial list. The same as the continuous control case, to our best knowledge, in all the existing literature treating optimal impulse control problems by dynamic programming principle, the terminal state is constraint free. A natural question aries: What if the terminal state is required to be constrained? Then one expects that, in general, the value function is not continuous, and might even not be defined somewhere. Recall that for continuous control case, there is a STLC condition that ensures the continuity of the value function. The major contribution of this paper is the discovery of an intrinsic condition that can play a similar role as STLC in the optimal impulse control problems. Under such a condition, the continuity of the value function will be proved and, consequently, it will be proved by mean of dynamic programming method that the value function is a viscosity solution to the Hamiton-Jacobi-Bellman QVI.

\ms

As far as applications are concerned, it is known that optimal impulse controls can be used in may areas, for examples, \cite{Bensoussan-Tapiero 1982} for management problems, \cite{Piunovskiy-Plakhov-Tumanov 2020} for SIR epidemic problems, \cite{Yadav-Balakrishnan 2006} for HIV treatment, \cite{Korn 1999} for mathematical finance, \cite{Leander-Lenhart-Protopopescu 2015} for some biology systems, to mention a few. For the optimal impulse control problem with a terminal state constraint, our motivation is as follow: Suppose a unit (could be a company, a bank, a state, or even a country) is running its business during a certain time period, say a month, a quarter-year, one year, etc. Besides it keeps normal running, at the end of the period, certain types of goods/assets (such as cash, food, gas, medicine, etc.) have to reach a certain reserve level. If it could not achieve the goal by its own production, it is allowed and has to buy from outside with some costs. This means that the state (goods/assets) of the unit can be  controlled by some impulses. The problem is to minimize the total cost with the terminal state constraint being satisfied. Clearly, such a framework is very general and could cover many real application problems. This also shows that the problem that we are going to study in the current paper is quite meaningful, both in mathematics and in applications.

\ms

Due to the presence of the terminal state constraint, the value function of the optimal impulse control problem is proved to be locally H\"older (or Lipschitz) continuous only, and it could grow at least linearly (no slower than the growth of the impulse cost). These properties essentially prevent us from directly using the current available techniques to prove the value function being the uniqueness of viscosity solution to the corresponding HJB QVI. Actually, we remind ourselves that the study of uniqueness of viscosity solutions to HJB equations for terminal state constrained problems was not successful in some other situations (see \cite{Soner-Touzi 2002, Touzi 2013} for the so-called stochastic target problems).

\ms

On the other hand, the obtained HJB QVI for the value function of our optimal impulse control problem with terminal state constraint looks like that for an optimal impulse control problem {\it without} terminal state constraint. Then a natural question arises: Can we indirectly characterize our value function by via the problem without constraint? In fact, if we can show that our impulse control problem is equivalent to a problem without terminal state constraint, then our goal is achieved. Some discussions will be carried out and the answer is not definite at the moment because we do not yet have the uniqueness of the viscosity solution to HJB QVI in the function class that our value function belongs to. Combining the above, we see that the issue of unique viscosity solution characterization for the value function of the problem with terminal state constraint remains challengingly open at the moment.

\ms

The rest of the paper is organized as follows. Section 2 is devoted to some preliminary results, including the formulation of the problem and the domain of the value function. In Section 3, we will investigate the continuity of the value function and present an interesting example there. In Section 4, we will derive dynamic programming principle and HJB QVI to which the value function is a viscosity solution. A comparison is made in Section 5 between the optimal impulse control problems with and without terminal state constraint, which reveals some interesting facts. Finally, concluding remarks are collected in Section 6.

\section{Preliminary Results}

In this section, we present some preliminary results.

\subsection{Formulation of the problem}

Let us first formulate our impulse control problem. Consider the following equation:
\bel{state1}X(s)=x+\int_t^sf(r,X(r))dr+\xi(s),\qq s\in[t,T],\ee
where $f:[0,T]\times\dbR^n\to\dbR^n$ is a given map, $(t,x)\in[0,T)\times\dbR^n$ is called an {\it initial pair}, and
\bel{xi}\xi(s)=\sum_{k\ges1}\xi_k{\bf1}_{[\t_k,T]}(s),\qq s\in[t,T]\ee
is called an {\it impulse control} with $\{\t_k\}_{k\ges1}\subset[t,T]$ being a non-decreasing finite sequence, and $\xi_k\in K$, $k\ges1$, called {\it admissible impulses}, for some non-empty closed convex cone $K\subseteq\dbR^n$ with the vertex at the origin. In the above, we allow $\t_k=\t_{k+1}$ for some $k\ges1$. Let $\sK[t,T]$ be the set of all impulse controls of form \rf{xi}. Under some mild conditions, for any initial pair $(t,x)\in[0,T)\times\dbR^n$ and impulse control $\xi(\cd)\in\sK[t,T]$, equation \rf{state1} admits a unique solution $X(\cd)\equiv X(\cd\,;t,x,\xi(\cd))$. Clearly, both $\xi(\cd)$ and $X(\cd)$ are right-continuous. In addition, we require that the terminal state satisfies the following constraint:
\bel{X(T)inD}X(T)\in\bar D,\ee
where $D$ is a non-empty proper domain in $\dbR^n$ (non-empty open and connected subset $D\ne\dbR^n$) with $\bar D$ being its closure. We may also call $\bar D$ a {\it target}. For any initial pair $(t,x)\in[0,T)\times\dbR^n$, we introduce the following associated {\it admissible impulse control} set
\bel{K_ad[t,T]}\sK^x[t,T]=\Big\{\xi(\cd)\in\sK[t,T]\bigm|X(T;t,x,\xi(\cd))\in \bar D\Big\}.\ee
In general, $\sK^x[t,T]$ is different from $\sK[t,T]$, and $\sK^x[t,T]$ could even be empty for some $(t,x)$. In the case $\sK^x[t,T]\ne\varnothing$, to measure the performance of the impulse control $\xi(\cd)$, we introduce the following cost functional
\bel{cost1}J(t,x;\xi(\cd))=\int_t^Tg(s,X(s))ds+h(X(T))+\sum_{k\ges1}\ell\big(\t_k,\wt X(\t_k-0),\xi_k\big),\ee
where
\bel{ghl}g:[0,T]\times\dbR^n\to[0,\infty),\q h:\dbR^n\to[0,\infty),\q\ell:[0,T]\times\dbR^n\times K\to(0,\infty)\ee
are suitable maps. Here, the terms on the right-hand side of \rf{cost1} are called the {\it running cost}, the {\it terminal cost} and the {\it impulse cost}, respectively. The meaning of $\wt X(\t_k-0)$ stands for the following: Suppose
$$\t_i<\t_{i+1}=\t_{i+2}=\cds=\t_{i'}<\t_{i'+1},$$
then
$$\wt X(\t_k-0)=X(\t_{i+1}-0)+\sum_{j=i+1}^{k-1}\xi_j,\qq i+1\les k\les i',\q\sum_{j=i+1}^i\xi_j\deq0,$$
which is the state right before the impulse $\xi_k$ is made. This is needed only if there are more than one separate impulses made at a same time (although such an impulse cannot be optimal). In the above, we may assume that $g$ and $h$ are just bounded uniformly from below. By a possible translation, we can simply assume that they are non-negative, for convenience. This will be assumed throughout of the paper. We emphasize that the impulse cost $\ell(t,x,\xi)$ is strictly positive. Mimicking the classical case, we formulate the following optimal impulse control problem.

\ms

{\bf Problem (C).} \rm For any initial pair $(t,x)\in[0,T)\times\dbR^n$, find a $\bar\xi(\cd)\in\sK^x[t,T]$ such that
\bel{infJ}J(t,x;\bar\xi(\cd))=\inf_{\xi(\cd)\in\sK^x[t,T]}J(t,x;\xi(\cd))=V(t,x).\ee

\ms

We call $\bar\xi(\cd)$ an {\it optimal impulse control}, the corresponding
$\bar X(\cd)\equiv X(\cd\,;t,x,\bar\xi(\cd))$ an {\it optimal state trajectory}, $(\bar X(\cd),\bar\xi(\cd))$ an {\it optimal pair}, and $V(\cd\,,\cd)$ the {\it value function} of Problem (C).

\ms

Recall a common convention that $\inf\varnothing=\infty$, regarding $\varnothing\subset\dbR$. Thus, it is convenient to make the following convention:
\bel{convention}J(t,x;\xi(\cd))=\infty,\qq\forall\xi(\cd)\in\sK[t,T]\setminus\sK^x[t,T].\ee
We let
\bel{D(V)}\sD(V)\equiv\Big\{(t,x)\in[0,T]\times\dbR^n\bigm|V(t,x)\hb{ is finite}\Big\},\ee
which is called the {\it domain} of the value function $V(\cd\,,\cd)$. Since $g(\cd\,,\cd),h(\cd)$ are non-negative and $\ell(\cd\,,\cd\,,\cd)$ is positive (see \rf{ghl}), one automatically has
\bel{D(V)in}\sD(V)=\Big\{(t,x)\in[0,T]\times\dbR^n\bigm|\sK^x[t,T]\ne
\varnothing\Big\}\equiv\sD(K;\bar D).\ee
The notation $\sD(K;\bar D)$ emphasizes the compatibility of the set $K$, $D$, and the dynamics \rf{state1}.

\subsection{Domain of the value function}

Before going further, let us first introduce the following hypotheses.

\ms

{\bf(H1)} $K\subseteq\dbR^n$ is a closed convex cone with the vertex located at the origin, and $D\subset\dbR^n$ is a non-empty proper convex domain (open and connected subset, different from $\dbR^n$).

\ms

{\bf(H2)} The map $f:[0,T]\times\dbR^n\to\dbR^n$ is continuous and there exists a constant $L>0$ such that
\bel{|f-f|}|f(t,x)-f(t,x')|\les L|x-x'|,\qq\forall t\in[0,T],~x,x'\in\dbR^n.\ee
\bel{|f|}|f(t,0)|\les L,\qq\forall t\in[0,T].\ee

\ms

Note that $K$ being a closed convex cone with the vertex located at the origin implies that if $\xi$ and $\xi'$ are two admissible impulses, so is $\xi+\xi'$. Also, $K$ might have empty interior. In what follows, we call the impulse control that contains no impulses the {\it trivial impulse control}, denote it by $\xi_0(\cd)$. Note that due to the presence of the (strictly positive) impulse cost, the trivial impulse control is different from the zero impulse control (which contains impulses with $\xi_k=0$). Let us first present the following result concerning the state trajectories.

\bp{state X} \sl Let {\rm(H1)--(H2)} hold. Then for any $(t,x)\in[0,T]\times\dbR^n$ and $\xi(\cd)\in\sK[t,T]$ of form \rf{xi}, state equation \rf{state1} admits a unique solution $X(\cd)=X(\cd\,;t,x,\xi(\cd))$, and the following estimates hold:
\bel{|X(s)|}|X(s)|\les e^{L(s-t)}(1+|x|)+\sum_{k\ges1}e^{L(s-\t_k)}|\xi_k|{\bf1}_{[\t_k,T]}(s),\qq s\in[t,T],\ee
\bel{|X(s')-X(s)|}|X(s'+0)-X(s+0)|\1n\les\1n L(s'\1n-\1n s)\1n+\1n\[e^{-Lt}(1\1n+\1n|x|)+\3n\sum_{\t_k\les s'}\2n|\xi_k|e^{-L\t_k}\]\big(e^{Ls'}\2n
-e^{Ls}\big)+\3n\2n\sum_{s<\t_k\les s'}\3n\2n|\xi_k|,\q0\les s<s'\les T.\ee
If $\h X(\cd)=X(\cd\,;t,\hat x,\xi(\cd))$ with $\hat x\in\dbR^n$, then
\bel{|X-X|}|X(s)-\h X(s)|\les e^{L(s-t)}|x-\hat x|,\qq s\in[t,T].\ee
\ep

\it Proof. \rm First of all, for any $(t,x)\in[0,T]\times\dbR^n$ and $\xi(\cd)\in\sK[t,T]$, by a standard argument, \rf{state1} admits a unique solution $X(\cd)=X(\cd\,;t,x,\xi(\cd))$. By Gronwall's inequality, we can get \rf{|X-X|}. We now prove \rf{|X(s)|}. From \rf{state1}, one has
$$|X(s)|\les|x|+\int_t^sL\big(1+|X(r)|\big)dr+|\xi(s)|\les|x|+L(s-t)+|\xi(s)|+L\int_t^s|X(r)|dr.$$
This implies
$$\ba{ll}
\ns\ds{d\over ds}\[e^{-L(s-t)}\int_t^s|X(r)|dr\]=e^{-L(s-t)}|X(s)|-e^{-L(s-t)}L\int_t^s|X(r)|dr\\
\ns\ds\qq\qq\qq\qq\qq\q~\les e^{-L(s-t)}\(|x|+L(s-t)+|\xi(s)|\).\ea$$
Hence,
$$e^{-L(s-t)}\int_t^s|X(r)|dr\les\int_t^se^{-L(r-t)}\(|x|+L(r-t)+|\xi(r)|\)dr.$$
Consequently,
$$\ba{ll}
\ns\ds|X(s)\les|x|+L(s-t)+|\xi(s)|+L\int_t^se^{L(s-r)}\(|x|+L(r-t)
+|\xi(r)|\)dr\\
\ns\ds\qq~=e^{L(s-t)}|x|+e^{L(s-t)}-1+|\xi(s)|+L\int_t^se^{L(s-r)}|\xi(r)|dr\\
\ns\ds\qq~\les e^{L(s-t)}(1+|x|)+\sum_{k\ges1}e^{L(s-\t_k)}|\xi_k|{\bf1}_{[\t_k,T]}(s).\ea$$
This proves \rf{|X(s)|}.

\ms

Now, let $t\les s<s'\les T$. Then
$$\ba{ll}
\ns\ds|X(s'+0)-X(s+0)|=\Big|\int_s^{s'}f(r,X(r))dr+\sum_{s<\t_k}\xi_k{\bf1}_{[\t_k,T]}(s')\Big|\\
\ns\ds\les\int_s^{s'}L\big(1+|X(r)|\big)dr+\sum_{s<\t_k}|\xi_k|{\bf1}_{[\t_k,T]}(s')\\
\ns\ds\les\int_s^{s'}L\(1+e^{L(r-t)}(1+|x|)+\sum_{k\ges1}e^{L(r-\t_k)}|\xi_k|
{\bf1}_{[\t_k,T]}(r)\)dr+\sum_{s<\t_k}|\xi_k|{\bf1}_{[\t_k,T]}(s')\\
\ns\ds\les L(s'-s)+\(e^{-Lt}(1+|x|)+\sum_{\t_k\les s'}|\xi_k|e^{-L\t_k}\)\big(e^{Ls'}
-e^{Ls}\big)+\sum_{s<\t_k\les s'}|\xi_k|.\ea$$
This proves \rf{|X(s')-X(s)|}. \endpf

\ms

In the above, \rf{|X(s)|} and \rf{|X-X|} are standard; whereas, \rf{|X(s')-X(s)|} seems to be new, from which, we see that although $s\mapsto X(s)$ might have jumps, these jumps can be controlled in some specific way. This is pretty natural and will be useful in the sequel.

\ms

To look at the domain $\sD(V)$ of the value function, we first note that under (H1), one always has
\bel{D(V)}\{T\}\times\bar D\subseteq\sD(V)=\sD(K;\bar D)\subseteq\sD(\dbR^n;\bar D)=[0,T]\times\dbR^n.\ee
Thus, $\sD(V)$ is always non-empty, and when $K=\dbR^n$, the domain $\sD(V)$ of $V(\cd\,,\cd)$ is the whole space $[0,T]\times\dbR^n$. Let us now present the following result.

\bp{Prop-3.2} \sl Let {\rm(H1)--(H2)} hold. Let
\bel{pa D in D-K}\pa D\subseteq D-K\equiv\{\eta-\xi\bigm|\eta\in D,~\xi\in K\},\ee
Then $\sD(V)$ is a non-empty open set in $[0,T]\times\dbR^n$.

\ep

\it Proof. \rm We already know that $\sD(V)$ is non-empty. Let $(t,x)\in\sD(V)$, then there exists an impulse control $\xi(\cd)\in\sK^x[t,T]$ such that $X(T;t,x,\xi(\cd))\in\bar D$. There are two cases.

\ms

\it Case 1. $X(T;t,x,\xi(\cd))\in D$. \rm Then there exists an $\e>0$ such that
$$B_\e\big(X(T;t,x,\xi(\cd))\big)\subseteq D,$$
where $B_\e(x)$ is the open ball centered at $x$ with radius $\e$.
Consequently, for $0<\d<e^{-L(T-t)}\e$, as long as $|x-\hat x|<\d$, one has
$$|X(T;t,\hat x,\xi(\cd))-X(T;t,x,\xi(\cd))|\les e^{L(T-t)}|x-\hat x|<\e.$$
Hence, $\xi(\cd)\in\sK^{\hat x}[t,T]$, leading to $(t,\hat x)\in\sD(V)$.  On the other hand, for $\hat t>t$, we let
$$\h\xi(\cd)=\sum_{\t_k\les\hat t}\xi_k{\bf1}_{[\hat t,T]}(\cd)+\sum_{\t_k>\hat t}\xi_k{\bf1}_{[\t_k,T]}(\cd).$$
This amounts to moving all the impulses no later than $\hat t$ to $\hat t$. Denote $\h X(\cd)=X(\cd\,;\hat t,x,\h\xi(\cd))$. Then for $s\in[\hat t,T]$,
$$\ba{ll}
\ns\ds|X(s)-\h X(s)|\les e^{L(s-\hat t)}\Big|X(\hat t)-x-\sum_{\t_k\les\hat t}\xi_k\Big|\\
\ns\ds\les e^{L(s-\hat t)}\int_t^{\hat t}|f(r,X(r))|dr\les Le^{L(s-\hat t)}\int_t^{\hat t}(1+|X(r)|)dr\\
\ns\ds\les Le^{L(s-\hat t)}\int_t^{\hat t}\(1+e^{L(r-t)}(1+|x|)+\sum_{t\les\t_k\les
r}e^{L(r-\t_k)}|\xi_k|\)dr\les C\(1+|x|+\sum_{t\les\t_k\les\hat t}|\xi_k|\)(\hat t-t).\ea$$
Hereafter, $C>0$ stands for a generic constant which could be different from line to line. Thus, when $\hat t-t>0$ small enough, we have $X(T;\hat t,x,\h\xi(\cd))\in D$, leading to
$\h\xi(\cd)\in\sK^x[\hat t,T]$. Finally, for $\hat t<t$, we take
$$\h\xi(\cd)=\sum_{k\ges1}\xi_k{\bf1}_{[\t_k,T]}(\cd),$$
i.e., we make a trivial extension of $\xi(\cd)$ from $[t,T]$ to $[\hat t,T]$. Denote
$\h X(\cd)=X(\cd\,;\hat t,x,\h\xi(\cd))$. Then
$$\ba{ll}
\ns\ds|\h X(t)-x|=|X(t;\hat t,x,\xi_0(\cd))-x|\les\int_{\hat t}^t|f(r,\h X(r))|dr
\les L\int_{\hat t}^t\(1+|\h X(r)|\)dr\\
\ns\ds\les L\int_{\hat t}^te^{L(r-\hat t)}(1+|x|)dr\les C(1+|x|)(t-\hat t).\ea$$
Hence, for $s\in[t,T]$,
$$|X(s)-\h X(s)|\les e^{L(s-t)}|x-\h X(t)|\les C(1+|x|)(t-\hat t).$$
Consequently, when $t-\hat t>0$ small enough, $X(T;\hat t,x,\h\xi(\cd))\in D$, leading to
$\h\xi(\cd)\in\sK^x[\hat t,T]$. Combining the above, we obtain
$$(\hat t,\h x)\in\sD(V),\qq\hb{if $|\hat t-t|+|\h x-x|$ is small enough}.$$

\ms

\it Case 2. $X(T;t,x,\xi(\cd))\in\pa D$. \rm In this case, by \rf{pa D in D-K}, there exists a $\bar\xi\in K$ such that by defining
$$\h\xi(\cd)=\sum_{k\ges1}\xi_k{\bf1}_{[\t_k,T]}(\cd)+\bar\xi{\bf1}_{\{T\}}(\cd),$$
we have
$$X(T;t,x,\h\xi(\cd))\in D.$$
Then it is reduced to Case 1. \endpf

\ms

The following result tells us more about $\sD(V)$.

\bp{Prop-2.2} \sl Let {\rm(H1)--(H2)} hold.

\ms

{\rm(i)} It holds that
\bel{barD-K=R}\bar D-K=\dbR^n,\ee
if and only if
\bel{D-K=R}D-K=\dbR^n.\ee
In this case,
\bel{sD=R}\sD(V)=\sD(K;\bar D)=[0,T]\times\dbR^n.\ee

{\rm(ii)} If $D$ is bounded, then \rf{barD-K=R} holds if and only if $K=\dbR^n$.

\ep

\it Proof. \rm (i) First of all, it is always true that $D-K\subseteq\bar D-K$. Thus, the sufficiency is clear. We now prove the necessity. Under (H1), both $\bar D-K$ and $D-K$ are convex. Moreover, $D-K$ is open and
\bel{D-K in D-K}\bar D-K\subseteq\cl{D-K}.\ee
In fact, for any $x_0\in D-K$, we have some $\eta_0\in D$ and $\xi_0\in K$ such that
$$x_0=\eta_0-\xi_0.$$
Since $D$ is open, there exists a $\d>0$ such that
$$\sO_\d(\eta_0)\equiv\{\eta\in\dbR^n\bigm||\eta-\eta_0|<\d\}\subseteq D.$$
Now, for any $x\in\sO_\d(x_0)$, we have
$$\eta\equiv\eta_0+x-x_0\in\sO_\delta(\eta_0)\subseteq D,$$
which leads to
$$x=x_0+x-x_0=\eta_0+(x-x_0)-\xi_0\equiv\eta-\xi_0\in D-K.$$
Thus, $D-K$ is open. The convexity of $\bar D-K$ and $D-K$ is clear. Next, for any
$$x=\eta-\xi\in\bar D-K,\qq\eta\in\bar D,\q\xi\in K,$$
we can find a sequence $\eta_k\in D$ such that $\eta_k\to\eta$. Then
$$x=\lim_{k\to\i}(\eta_k-\xi)\in\cl{D-K},$$
proving \rf{D-K in D-K}. Now, if $D-K\ne\bar D-K$, by the convexity of $D-K$, we must have $\cl{D-K}\ne\dbR^n$. This, together with \rf{D-K in D-K}, contradicts \rf{barD-K=R}.

\ms

Finally, we prove \rf{sD=R}. For any $(t,x)\in[0,T]\times\dbR^n$, under the trivial impulse control $\xi_0(\cd)$, the state will arrive at $X(T-0;t,x,\xi_0(\cd))\in\dbR^n$. By \rf{barD-K=R}, we have some $\eta\in\bar D$ and $\xi\in K$ such that
$$X(T-0;t,x,\xi_0(\cd))=\eta-\xi.$$
Then by defining impulse control
$$\h\xi(\cd)=\xi{\bf1}_{\{T\}}(\cd),$$
we have
$$X(T;t,x,\h\xi(\cd))=X(T-0;t,x,\xi_0(\cd))+\xi=\eta\in\bar D.$$
Thus, $\sK^x[t,T]\ne\varnothing$. This proves our conclusion (see \rf{D(V)in}).

\ms

(ii) First of all, if $K=\dbR^n$, then of course. \rf{barD-K=R} holds. Now, if $K\ne\dbR^n$, then there must be a $\z\in K$, $|\z|=1$ such that
$$\lan\z,\xi\ran\ges0,\qq\forall\xi\in K.$$
Now, we claim that $\l\z\notin\bar D-K$ for large enough $\l>0$. In fact, if there exists an $\eta^\l\in\bar D$ and a $\xi^\l\in K$ such that
$$\l\z=\eta^\l-\xi^\l.$$
This leads to
$$\eta^\l=\l\z+\xi^\l.$$
Hence,
$$\lan\eta^\l,\z\ran=\l+\lan\xi^\l,\z\ran\ges\l.$$
Since $\{\eta^\l\}_{\l>0}$ is bounded, we may assume that $\eta^\l\to\bar\eta$. But this will lead to a contradiction. \endpf

\ms

One of the most interesting examples satisfying \rf{barD-K=R} is the following:
$$\bar D=K=\dbR^n_+\equiv\big\{x\in\dbR^n\bigm|x_i\ges0\big\}.$$
Consequently, for such a case, one has \rf{sD=R}. The above proposition gives two important cases: $D$ is a bounded set with $K=\dbR^n$ and $D$ is unbounded with $K\ne\dbR^n$ such that \eqref{barD-K=R} holds. They are  mutually exclusive. However, we point out that they are not exhausting. Here is a simple example that is neither of the above cases: In $\dbR^2$, let
$$K=\{(x_1,x_2)\bigm|x_1,x_2\ges0\},\qq D=\Big\{(x_1,x_2)\bigm|x_2>{1\over1-x_1},~x_1<1\Big\}.$$
Then $D$ is unbounded and
$$\bar D\subseteq\bar D-K=\big\{(x_1,x_2)\bigm|x_1<1\big\}\ne\dbR^2.$$
If we regard \rf{barD-K=R} as the best case, then since $0\in K$, the worst case should be
\bel{D-K=D}\bar D-K=\bar D.\ee
An example of such is the following:
\bel{KD}K=\{(x_1,x_2)\bigm|x_1,x_2\ges0\},\qq D=\{(x_1,x_2)\bigm|x_1,x_2<0\}.\ee
For the case that \rf{barD-K=R} fails, including the case of \rf{D-K=D}, when $X(T-0)\notin\bar D$, there is no way to make an impulse at $T$ so that $X(T)\in\bar D$. Therefore,
$$\big[\big\{T\}\times(\bar D-K)^c\big]\cap\sD(V)=\varnothing,\q\hb{or}\q\big[\{T\}\times\dbR^n\big]\cap\sD(V)
\ne\{T\}\times\dbR^n.$$
For such a case, one has to make impulses before $T$ and drive the state to $\bar D$ (at $T$) via the state equation. Thus, it might still be possible that
$$\big[\{t\}\times\dbR^n\big]\cap\sD(V)=\{t\}\times\dbR^n,$$
for some $t\in[0,T)$. We will see a concrete example below.

\ms

Let us now look at the following simple example to get some more feeling.

\bex{Ex-2.4} \rm Consider
\bel{2.13}X(s)=x+(s-t)+\xi(s),\qq s\in[t,T].\ee
We consider several cases.

\ms

(i) $K=[0,\i)$, $D=(0,1)$. For this case,
$$\bar D-K=[0,1]-[0,\i)=(-\i,1],$$
and
$$\sD(V)=\big\{(t,x)\in[0,T]\times\dbR\bigm|x+T-t\les1\big\}.$$

\ms

(ii) $K=(-\i,0]$, $D=(0,1)$. For this case,
$$\bar D-K=[0,1]-(-\i,0]=[0,\i),$$
and
$$\sD(V)=\big\{(t,x)\in[0,T]\times\dbR\bigm|x+T-t\ges0\big\}.$$

(iii) $K=[0,\i)$, $D=(0,\i)$. For this case,
$$\bar D-K=[0,\i)-[0,\i)=\dbR,\qq\sD(V)=[0,T]\times\dbR.$$

(iv) $K=[0,\i)$, $D=(-\i,0)$. For this case,
$$\bar D-K=(-\i,0]-[0,\i)=(-\i,0]=\bar D,$$
and
$$\sD(V)=\big\{(t,x)\in[0,T]\times\dbR\bigm|x+T-t\les0\big\}.$$

\ex

Now, we look at the general situation. Under (H1)--(H2), for any $x\in D$, we may let $B_\e(x)\subseteq D$, with $B_\e(x)$ being the ball centered at $x$ with radius $\e$. Let $t\in[0,T]$ such that $T-t>0$ is small enough so that
$$|X(s;t,x,\xi_0(\cd))-x|\les(1+L)e^{LT}(1+|x|)(T-t)<\e,\qq s\in[t,T].$$
This means that $\xi_0(\cd)\in\cK^x[t,T]$. Hence, under (H1)--(H2), the following is always true:
\bel{non-empty}\sD(V)=\sD(K;\bar D)\ne\varnothing.\ee
We now would like to get a more precise description of $\sD(V)$. For state equation \rf{state1}, we consider the following ``backward'' system
$$Y(s)=\z-\int_s^Tf(r,Y(r))dr,\qq s\in[0,T],$$
with $\z\in\bar D-K$. The solution is denoted by $Y(\cd\,;T,\z)$. Let
$$\sY(t;T,\bar D-K)=\{Y(t;T,\z)\bigm|\z\in\bar D-K\}.$$
For any $x\in\sY(t;T,\bar D-K)$, one has some $\z\in\bar D-K$ such that
$$x=Y(t;T,\z).$$
Then, with the trivial impulse control $\xi_0(\cd)$, we have
$$X(T-0;t,x,\xi_0(\cd))=Y(T;T,\z)=\z\in\bar D-K.$$
One can choose some $\xi\in K$ such that
$$X(T-0;t,x,\xi_0(\cd))+\xi\in\bar D.$$
Thus, $\sY(t;T,\bar D-K)$ is the set of all possible initial state that if the system starts at $(t,x)$, the state will reach $\bar D-K$ at $T$ under $\xi_0(\cd)$. Then, under a possible impulse at $T$, the state will hit $\bar D$. Now, let $\Pi=\{t_0,t_1,\cds,t_N\}$ be a partition of $[t,T]$ with $t=t_0<t_1<t_2<\cds<t_N=T$. Then we may define inductively
$$\ba{ll}
\ns\ds\sY^N_\Pi\equiv\bar D,\\
\ns\ds\sY^{N-1}_\Pi=\sY(t_{N-1};t_N,\bar D-K)\equiv\sY(t_{N-1};t_N,\sY^N_\Pi-K),\\
\ns\ds\sY^{N-2}_\Pi=\sY(t_{N-2};t_{N-1},\sY^{N-1}_\Pi-K),\\
\ns\ds\cds\cds\\
\ns\ds\sY^1_\Pi=\sY(t_1;t_2\sY^2_\Pi-K),\\
\ns\ds\sY^0_\Pi=\sY(t_0;t_1\sY^1_\Pi-K).\ea$$
We denote
$$\sY(t;\Pi)=\sY^0_\Pi-K,$$
which is the set of all initial states $x\in\dbR^n$ such that if the system starts at $(t,x)$, with possible impulses at $t_0,t_1,\cds,t_N$, the state will reach $\bar D$ at $T$. This can be described by the following:
$$\ba{ll}
\ns\ds\qq\qq\qq\qq\bar D\equiv\sY^N_\Pi~{\color{red}\La}~\sY^N_\Pi-K\\
\ns\ds\qq\qq\qq\qq\qq\qq\qq\qq{\color{blue}\Uparrow}\\
\ns\ds\qq\qq\sY^{N-1}_\Pi-K~{\color{red}\Ra}~\sY^{N-1}_\Pi\equiv\sY(t_{N-1};t_N,\sY^N_\Pi-K)\\
\ns\ds\qq\qq\qq{\color{blue}\Uparrow}\\
\ns\ds\sY(t_{N-2};t_{N-1},\sY^{N-1}_\Pi-K)\equiv\sY^{N-2}_\Pi~{\color{red}\La}~
\sY^{N-2}_\Pi-K\\
\ns\ds\qq\qq\qq\cds\cds\cds\cds\cds\\
\ns\ds\qq\qq\qq\qq\qq\qq\qq\qq\qq{\color{blue}\Uparrow}\\
\ns\ds\qq\qq\sY^1_\Pi-K~{\color{red}\Ra}~\sY^1_\Pi\equiv\sY(t_1;t_2,\sY^2_\Pi-K)\\
\ns\ds\qq\qq\qq{\color{blue}\Uparrow}\\
\ns\ds\sY(t_0;t_1,\sY^1_\Pi-K)\equiv\sY^0_\Pi~{\color{red}\La}~
\sY^0_\Pi-K\ea$$
In the above, horizontal arrows represent making impulses, and upper arrows represent running state equations. Clearly, for any two partitions $\Pi_1$ and $\Pi_2$ of $[t,T]$ with $\Pi_1\subseteq\Pi_2$, i.e., $\Pi_2$ is a refinement of $\Pi_1$, we have
$$\sY(t;\Pi_1)\subseteq\sY(t;\Pi_2).$$
Hence, we may define
$$\sY(t)=\bigcup_{\Pi}\sY(t;\Pi)=\lim_{\|\Pi\|\to0}\sY(t;\Pi),$$
where $\|\Pi\|$ is the mesh size of $\Pi$ defined by
$$\|\Pi\|=\max_{1\les k\les N}(t_k-t_{k-1}).$$
From the construction, we see that $\sY(t)$ is the set of all initial state that if the system starts from $(t,x)$, then with impulse controls, the state can reach $\bar D$ at $T$, i.e.,
\bel{}\sK^x[t,T]\ne\varnothing\qq\iff\qq x\in\sY(t).\ee
Hence, we have the following characterization of $\sD(V)$:
\bel{D(V)*}\cl{\sD(V)}=\cl{\bigcup_{t\in[0,T]}\big[\{t\}\times\sY(t)\big]}.\ee
The following example gives a concrete construction of $\sY(t)$.

\bex{2.4} \rm Consider
$$\left\{\2n\ba{ll}
\ns\ds X_1(s)=x_1+\int_t^sX_2(r)dr+\xi_1(s),\\
\ns\ds X_2(s)=x_2-\int_t^sX_1(r)dr+\xi_2(s),\ea\right.\qq s\in[t,T].$$
Let
$$D=\{(x_1,x_2)\bigm|x_1^2+x_2^2<1\},\qq K=\dbR^2_+\equiv\{(x_1,x_2)\bigm|x_1,x_2\ges0\}.$$
The backward system reads
$$\left\{\2n\ba{ll}
\ns\ds Y_1(s)=\eta_1-\int_s^TY_2(r)dr-\xi_1(s),\\
\ns\ds Y_2(s)=\eta_2+\int_s^TY_1(r)dr-\xi_2(s),\ea\right.\qq s\in[0,T].$$
For any $\eta=(\eta_1,\eta_2)\in\bar D$ and $(\z_1,\z_2)\in K$, let
$$\begin{pmatrix}Y_1(s)\\ Y_2(s)\end{pmatrix}=\begin{pmatrix}\cos(T-s)&-\sin(T-s)\\\sin(T-s)&\cos(T-s)\end{pmatrix}
\begin{pmatrix}\eta_1-\z_1\\ \eta_2-\z_2\end{pmatrix}.$$
Note that as $s$ decreases from $T$, the vector $(Y_1(s),Y_2(s))^\top$ turns counter-clockwise. We may keep making impulses to see that
$$\ba{ll}
\ns\ds\sY(T)=\co\(\bar D\cup\{(x_1,x_2)\bigm|x_1,x_2\les0\}\),\\
\ns\ds\sY(t)=\co\(\sY(T)\cup\{(x_1,x_2)\bigm|x_2\les x_1\cot(T-t),~x_1>0\}\),\qq0\les T-t<{\pi\over4},\\
\ns\ds\sY(T-{\pi\over2})=\{(x_1,x_2)\bigm|x_2\les1\},\qq\sY(t)=\dbR^2,\qq T-t>{\pi\over2},\ea
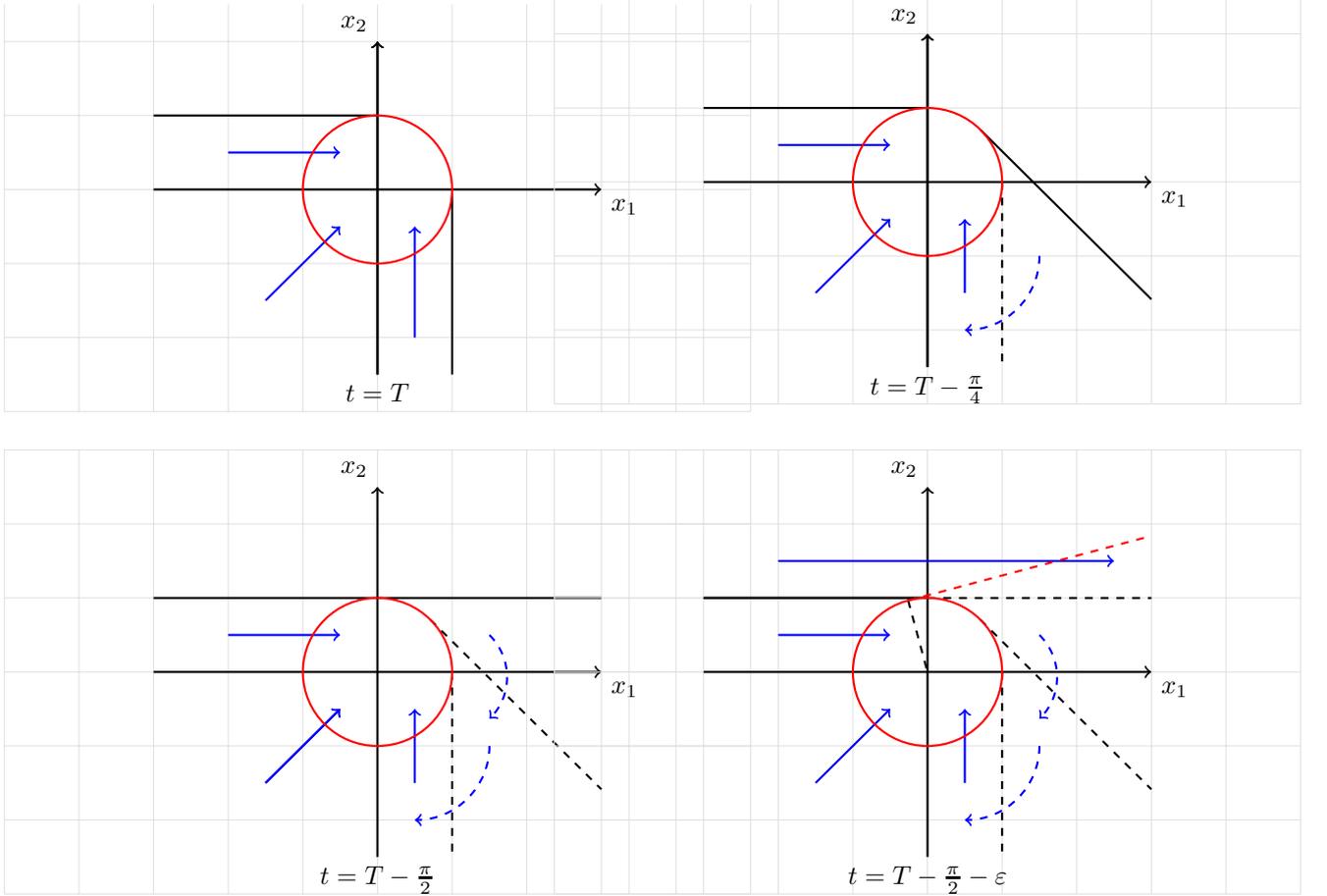
$$
$$
\begin{tikzpicture}
\draw[step=1cm,gray!25!,very thin] (-5,-3) grid (5,2.5);
\draw[thick,->] (-3,0) -- (3,0) node[anchor=north west] {$x_1$};
\draw[thick,->] (0,-2.5) -- (0,2) node[anchor=south east] {$x_2$};
\draw[blue,thick,->] (-2,0.5) -- (-0.5,0.5); 
\draw[blue,thick,->] (0.5,-2) -- (0.5,-0.5); 
\draw[thick,->] (0,-2.5) -- (0,2); 
\draw[thick] (1,0) -- (1,-2.5); 
\draw[thick] (-3,1) -- (0,1); 
\draw[blue,thick,->] (-1.5,-1.5) -- (-0.5,-0.5);
\draw[red, thick] (0,0) circle (1 cm);
\draw (0,-2.5) node[below]{$t=T$};
\end{tikzpicture}\5n\5n\5n\5n\5n\5n\5n\5n\5n
\begin{tikzpicture}
\draw[step=1cm,gray!25!,very thin] (-5,-3) grid (5,2.5);
\draw[thick,->] (-3,0) -- (3,0) node[anchor=north west] {$x_1$};
\draw[thick,->] (0,-2.5) -- (0,2) node[anchor=south east] {$x_2$};
\draw[blue,thick,->] (-2,0.5) -- (-0.5,0.5); 
\draw[blue,thick,->] (0.5,-1.5) -- (0.5,-0.5); 
\draw[thick,->] (0,-2.5) -- (0,2); 
\draw[dashed,thick] (1,0) -- (1,-2.5); 
\draw[thick] (-3,1) -- (0,1); 
\draw[thick] (0.7071,0.7071) -- (3,-1.5858);
\draw[blue,thick,->] (-1.5,-1.5) -- (-0.5,-0.5);
\draw[red, thick] (0,0) circle (1cm);
\draw[blue,dashed,thick, ->] (1.5,-1) arc (0:-90:1);
\draw (0,-2.5) node[below]{$t=T-{\pi\over4}$};
\end{tikzpicture}$$
\vskip-0.5cm
$$
\begin{tikzpicture}
\draw[step=1cm,gray!25!,very thin] (-5,-3) grid (5,3);
\draw[thick,->] (-3,0) -- (3,0) node[anchor=north west] {$x_1$};
\draw[thick,->] (0,-2.5) -- (0,2.5) node[anchor=south east] {$x_2$};
\draw[blue,thick,->] (-2,0.5) -- (-0.5,0.5); 
\draw[blue,thick,->] (0.5,-1.5) -- (0.5,-0.5); 
\draw[dashed,thick] (1,0) -- (1,-2.5); 
\draw[thick] (-3,1) -- (3,1); 
\draw[blue,thick,->] (-1.5,-1.5) -- (-0.5,-0.5);
\draw[dashed,thick] (0.7071,0.7071) -- (3,-1.5858);
\draw[blue,thick,->] (-1.5,-1.5) -- (-0.5,-0.5);
\draw[red, thick] (0,0) circle (1cm);
\draw[blue,dashed,thick, ->] (1.5,-1) arc (0:-90:1);
\draw[blue,dashed,thick, ->] (1.5,0.5) arc (45:-45:0.8);
\draw (0,-2.5) node[below]{$t=T-{\pi\over2}$};
\end{tikzpicture}\5n\5n\5n\5n\5n\5n\5n\5n\5n
\begin{tikzpicture}
\draw[step=1cm,gray!25!,very thin] (-5,-3) grid (5,3);
\draw[thick,->] (-3,0) -- (3,0) node[anchor=north west] {$x_1$};
\draw[thick,->] (0,-2.5) -- (0,2.5) node[anchor=south east] {$x_2$};
\draw[blue,thick,->] (-2,0.5) -- (-0.5,0.5); 
\draw[blue,thick,->] (0.5,-1.5) -- (0.5,-0.5); 
\draw[blue,thick,->] (-2,1.5) -- (2.5,1.5); 
\draw[dashed,thick] (1,0) -- (1,-2.5); 
\draw[thick] (-3,1) -- (0,1); 
\draw[dashed,thick] (0,1) -- (3,1);
\draw[dashed,thick] (0.7071,0.7071) -- (3,-1.5858);
\draw[thick] (-3,1) -- (0,1);
\draw[red,dashed,thick] (-0.259,0.966) -- (3,1.839);
\draw[dashed,thick] (-0.259,0.966) -- (0,0);
\draw[blue,thick,->] (-1.5,-1.5) -- (-0.5,-0.5);
\draw[red, thick] (0,0) circle (1cm);
\draw[blue,dashed,thick, ->] (1.5,-1) arc (0:-90:1);
\draw[blue,dashed,thick, ->] (1.5,0.5) arc (45:-45:0.8);
\draw (0,-2.5) node[below]{$t=T-{\pi\over2}-\e$};
\end{tikzpicture}$$
where $\co(M)$ is the convex hull of the set $M$, i.e., the smallest convex set containing $M$. In the illustrative figures, the blue arrow lines give the directions of impulses; the dashed arcs give the directions of the points turning. Thus, in the last figure (of the situation $T-t>{\pi\over2}+\e$, any initial point $(x_1,x_2)$ with $x_2>1$, one could first make a horizontal impulse $\xi=(\xi_1,0)$ so that $(x_1+\xi_1,x_2)$ is on the right of the dashed red line. Then by the original system which makes the point turning clockwise, and at $t=T$, the point will be in $\sY(T)$. By making an impulse at $t=T$, the state will get into $\bar D$.

\ms

From the above, we see that
$$\ba{ll}
\ns\ds\sD(V)\cap\big[\{t\}\times\dbR^n\big]=\sD(K;\bar D)\cap\big[\{t\}\times\dbR^n\big]\ne\{t\}\times\dbR^n,\qq0\les T-t\les{\pi\over2},\\
\ns\ds\sD(V)\cap\big[\{t\}\times\dbR^n\big]=\sD(K;\bar D)\cap\big[\{t\}\times\dbR^n\big]=\{t\}\times\dbR^n,\qq T-t>{\pi\over2}.\ea$$
However, one has
$$\sD(V)=\sD(K;\bar D)\supseteq[0,T]\times\bar D.$$
\ex

The above \rf{D(V)*} gives a characterization for the domain $\sD(V)$ of the value function, in some sense. More needs to be done. For example, what will be the boundary of such a domain, how the value function behaves near the boundary of this domain, and so on. We are going to leave these problems open for the time-being, and hope that we will be able to report some relevant results in the near future. Instead, in the current paper, we consider the things more relevant to the continuity of the value function. Let us make some more preparations.

\bl{Lemma 2.6} \sl Let {\rm(H1)} and \rf{D-K=R} hold. Then there exists a nondecreasing continuous function $\n:(0,\infty)\to(0,\infty)$ only depending on $K$ and $D$ such that for any $x\in\dbR^n$, there exists a $\xi\in K$ satisfying
\bel{x+xi in D}x+\xi\in\bar D,\qq|\xi|\les\n(|x|).\ee
Further, if $D$ is bounded (in this case, it is necessary that $K=\dbR^n$) or if $D$ is a conic domain with the vertex located at the origin, then for some constant $C_0$, only depending on $D$ and $K$,
\bel{m=1+x}\n(r)=C_0(1+r).\ee

\el

\it Proof. \rm First we claim that for any fixed $M>0$, there exists a $\n_M>0$ only depending on $K$ and $D$ such that for any $x\in\dbR^n$ with $|x|\les M$, there exists a $\xi\in K$ with $|\xi|\les\n_M$ satisfying $x+\xi\in\bar D$, i.e.,
$$d(x+\xi,\bar D)\equiv\inf_{\eta\in\bar D}|x+\xi-\eta|=0.$$
Suppose this claim fails. Then there exists some $M>0$ such that one can find a sequence $x_k\in\dbR^n$ with $|x_k|\les M$, and
$$\inf_{\eta\in\bar D}|x_k+\xi-\eta|=d(x_k+\xi,\bar D)\ges\d_k,\qq\forall\xi\in K,~|\xi|\les k,~k\ges1,$$
for some $\d_k>0$. We may assume $x_k\to x_0$. Now for $x_0$, by \rf{D-K=R}, there exists a $\xi_0\in K$ such that
$$\eta_0=x_0+\xi_0\in D.$$
Hence, for $k>|\xi_0|$, large enough, we have (noting that $D$ is open)
$$x_k+\xi_0=x_k-x_0+x_0+\xi_0=x_k-x_0+\eta_0\in D\subseteq\bar D.$$
This is a contradiction. Hence, the claim is true. Now, we define
$$\n(x)=\left\{\2n\ba{ll}
\ds\n_{k-1}+2(\n_k-\n_{k-1})(|x|-k+1),\qq&k-1\les|x|\les k-{1\over2},\\
\ns\ds\n_k,\qq&k-{1\over2}\les|x|\les k,\ea\right.\qq k\ges1.$$
which satisfies the our requirement.

\ms

Next, if $D$ is bounded with the bound $M>0$, then for any $x\in\dbR^n$, let $\xi\in K=\dbR^n$ such that $x+\xi=\eta\in D$. Clearly,
$$|\xi|\les|x|+|\eta|\les M+|x|,$$
which is what we want.

\ms

Finally, let $D$ be a conic domain with vertex located at the origin. Then
under our condition \rf{D-K=R}, we can find a $\xi_0\in D\cap K$, $\xi_0\ne0$. (Note here that $0\notin D$.) Now, for any $x\in\dbR^n$, if $x\in\bar D$ (which includes the case that $x=\l\xi_0$ for some $\l\ges0$), we trivially have $x=x-0\in\bar D-K$. Hence, by taking $\eta=x$ and $\xi=0$, we have \rf{x+xi in D}--\rf{m=1+x}. Next, let $x\notin \bar D$. If
$$x=-\l\xi_0,$$
for some $\l>0$, then by taking $\eta=0\in\bar D$ and $\xi=\l\xi_0$, we have \rf{x+xi in D}--\rf{m=1+x}.
Hence, we need only to look at the case that $x\notin\bar D$, with $x$ and $\xi_0$ being linearly independent. Consider the two-dimensional space $\dbH$ spanned by $x$ and $\xi_0$. After a proper linear transformation, we may assume the following situation in $\dbR^2$:
$$\ba{ll}
\ns\ds D\cap K\ni\xi_0=(1,0)^\top,\q\{-\l\xi_0\bigm|\l\ges0\}=\{(-\l,0)^\top\bigm|\l\ges0\}=(-\infty,0]
\times\{0\},\\
\ns\ds\bar D\cap\dbH=\{(\eta_1,\eta_2)^\top\bigm|-\a\eta_1\les\eta_2\les\b \eta_1,~\eta_1\ges0\big\},\ea$$
for some $\a,\b>0$ (depending on $D$ and $K$). Then for $(x_1,x_2)\in\dbR^2\setminus\bar D$, in the case that $x_2>0$, we must have $x_2>\b x_1$. Take
\bel{eta}\eta_1={x_2\over\b}>0,\q\eta_2=x_2,\q\l={x_2\over\b}-x_1>0.\ee
One has $\eta\equiv(\eta_1,\eta_2)^\top\in\bar D\cap\dbH\subseteq\bar D$, $\xi\equiv\l\xi_0\in K$, and
$$(\eta_1,\eta_2)-\l(1,0)=\({x_2\over\b},x_2\)-\({x_2\over\b}-x_1\)(1,0)
=(x_1,x_2)=x.$$
Clearly,
$$|\xi|=\l\les{\sqrt{1+\b^2}\over\b}|x|,\q|\eta|={\sqrt{1+\b^2}\over\b}|x_2|
\les{\sqrt{1+\b^2}\over\b}|x|.$$
Thus, \rf{x+xi in D}--\rf{m=1+x} hold. Likewise, in the case $x_2<0$, we must have $x_2<-\a x_1$. Take
\bel{eta*}\eta_1={-x_2\over\a}>0,\q\eta_2=x_2,\q\l=-{x_2\over\a}-x_1>0.\ee
Then $\eta=(\eta_1,\eta_2)^\top\in\bar D\cap\dbH\subseteq\bar D$, $\xi=\l\xi_0
\in K$, and
$$(\eta_1,\eta_2)-\l(1,0)=\({-x_2\over\a},x_2\)+\({x_2\over\a}+x_1\)(1,0)
=(x_1,x_2)=x.$$
Also,
$$|\xi|=\l\les{\sqrt{1+\a^2}\over\a}|x|,\q|\eta|={\sqrt{1+\a^2}\over\a}|x_2|\les
{\sqrt{1+\a^2}\over\a}|x|.$$
Hence, \rf{x+xi in D}--\rf{m=1+x} hold for this case as well.  \endpf

\ms

The point of the above result is that under condition \rf{D-K=R}, the size of the impulse $\xi$ that drives the state $x$ to the constraint $D$ is controlled by $|x|$. This is very important below. We also note that \rf{D-K=R} seems to be a little stronger than \rf{barD-K=R}. In the case that $K^\circ\ne\varnothing$ ($K^\circ$ is the interior of $K$), they are equivalent. It is not clear to us at the moment if they are equivalent in general. Also, we point out that the function $r\to\n(r)$ can grow arbitrarily fast. Here is a simple example.

\bex{} \rm Let $\n_0:[0,\i)\to[0,\i)$ be a continuous strictly increasing function with $\n_0(0)=0$ and $\n_0(r)\to\i$ as $r\to\i$ (for example, $\n_0(r)=e^r-1$). Let
$$D=\big\{(x_1,x_2)\in\dbR^2\bigm|x_1>\n_0(|x_2|)\big\},\qq K=\big\{(x_1.x_2)\in\dbR^2\bigm|x_1,x_2\ges0\big\}.$$
Then, for any $(x_1,x_2)\in\dbR^2\setminus\bar D$, one can take $\xi=(\xi_1,\xi_2)\in K$ with
$$\xi_1\ges\n_0(|x_2|)-x_1,\qq\xi_2=0,$$
which will lead to $x+\xi\in\bar D$. This is actually the best choice as far as the norm $|\xi|$ is concerned. Clearly, if $\n(\cd)$ is the function appeared in Lemma \ref{Lemma 2.6}, then with $x_1=0$, one has
$$\n_0(|x_2|)\les\xi_1\les|\xi|\les\n(|x_2|),\qq\forall x_2\in\dbR.$$
Thus, $\n(\cd)$ cannot be growing slower than $\n_0(\cd)$.

\ex

To conclude this section, let us present one more example.

\bex{2.8} \rm Let $K=\dbR^2$, and
$$D=\{(x_1,x_2)\in\dbR^2\bigm|x_1\ges0\},$$
then for any $x=(x_1,0)\in\dbR^2$ with $x_1<0$, we see that the $\xi=(\xi_1,\xi_2)\in K$ that makes $x+K\in\bar D$ with smallest possible $|\xi|$ should be $\xi=(-x_1,0)$ and $|\xi|=|x|$.

\ex

The point that we want to make in the above example is that as long as $D\ne\dbR^n$ and $D$ is convex, the function $\n(r)$ appears in Lemma \ref{Lemma 2.6} will be at least of linear growth.

\section{Properties of the Value Functions}

In this section, we will present some properties of the value function $V(\cd\,,\cd)$, including its continuity.

\subsection{Some bounds}

First, we introduce the following strengthened (H1).

\ms

{\bf(H1$'$)} Let $D\subset\dbR^n$ be a bounded domain and $K=\dbR^n$ or $D\subset\dbR^n$ be a conic domain with the vertex being at the origin and $K\subseteq\dbR^n$ be a closed convex cone with the vertex located at the origin such that \rf{D-K=R} holds.

\ms
Next, let us introduce the following additional hypotheses.

\ms

{\bf(H3)} Maps $g:[0,T]\times\dbR^n\to[0,\infty)$ and $h:\dbR^n\to[0,\infty)$ are continuous. There exist constants $L,\m>0$ and $0<\d\les1$ such that
\bel{g,h}0\les g(t,x),h(x)\les L\big(1+|x|^{\m+\d}\big),\qq\forall(t,x)\in[0,T]\times\dbR^n,\ee
\bel{g-g,h-h}|g(t,x)-g(t,x')|,|h(x)-h(x')|\les L\(1+|x|^\m\vee|x'|^\m\)|x-x'|^\d,\qq\forall t\in[0,T],x,x'\in\dbR^n.\ee

\ms

{\bf(H4)} Map $\ell:[0,T]\times\dbR^n\times K\to(0,\infty)$ is continuous. There exist constants $\ell_0,\d_0>0$, $\b\in(0,1]$ and $L,\m,\d>0$ the same as those in (H3), as well as a continuous decreasing functions $\a_0,\a:[0,T]\to(0,\infty)$ such that
\bel{ell(1)}\ell_0+\a_0(t)|\xi|^\b\les\ell(t,x,\xi)\les L+\a(t)|\xi|^\b,\qq(t,x,\xi)\in[0,T]\times\dbR^n\times K,\ee
\bel{ell(2)}|\ell(t,x,\xi)-\ell(t,x',\xi)|\les L\(1+|x|^\m\vee|x'|^\m\)|x-x'|^\d,\qq\forall(t,\xi)\in[0,T]\times
K,~x,x'\in\dbR^n,\ee
\bel{ell(3)}\ba{ll}
\ns\ds\min\big\{\ell(t,x,\xi)+\ell(t,x+\xi,\xi'),\ell(t,x,\xi')
+\ell(t,x+\xi',\xi)\big\}-\ell(t,x,\xi+\xi')\ges\d_0,\\
\ns\ds\qq\qq\qq\qq\qq(t,x)\in[0,T]\times\dbR^n,~\xi,\xi'\in K,\ea\ee
\bel{ell(4)}\ell(t,x,\xi)-L|t-t'|\les\ell(t',x,\xi)\les\ell(t,x,\xi),\qq t\les t',~x\in\dbR^n,~\xi\in K.\ee

As we have indicated in the introduction section, one can assume that $g$ and $h$ are bounded below uniformly. Here, we directly assume them to be non-negative just for convenience. Condition \rf{ell(1)} implies that as long as an impulse is made, no matter how small the $\xi$ is, there is a strictly positive fixed cost $\ell_0$. Also, roughly speaking, the larger the $|\xi|$, the larger the cost. Condition \rf{ell(3)} means that if at $(t,x)$ an impulse of size $\xi+\xi'$ needs to be made, then one should make just one impulse of that size instead of making an impulse of size $\xi$ immediately followed by another with size $\xi'$. Hence, in an optimal impulse control, $\t_k<\t_{k+1}$ if both are impulsive moments. In the case that $\ell(t,x,\xi)$ is independent of $x$, this condition is reduced to
$$\ell(t,\xi+\xi')<\ell(t,\xi)+\ell(t,\xi'),$$
which is a classical condition assumed in the optimal impulse control problems. Because of this condition, $\xi\mapsto\ell(t,x,\xi)$ should be ``sublinear''. Hence, $\b\in(0,1]$ and $\xi\mapsto\ell(t,x,\xi)$ grows at most linearly (see \rf{ell(1)}). Condition \rf{ell(4)} means that if an impulse is going to be made, then the later the better, which is essentially due to the discount effect.

\ms

Our goal in this section is to obtain, under certain conditions, including (H1)--(H4), the bounds of the value functions, the smaller class of impulse controls on which the value functions are the infimum of the cost functional, and each impulse control in this smaller class has no more than a fixed number of impulses with the sizes of the impulses being bounded. More precisely, we have the following result.

\bp{Prop-3.2} \sl Let {\rm(H1)--(H4)} and \rf{D-K=R} hold. Then $\sD(V)=[0,T]\times\dbR^n$ and
\bel{V<}0\les V(t,x)\les\bar\n(|x|),\qq\forall(t,x)\in
[0,T]\times\dbR^n,\ee
for some continuous increasing function $\bar\n:[0,\i)\to[0,\i)$.
Moreover, for any $(t,x)\in[0,T]\times\dbR^n$, Problem {\rm(C)} admits an optimal impulse control, and
\bel{V=}V(t,x)=\inf_{\xi(\cd)\in\sK^x_0[t,T]}J(t,x;\xi(\cd)),\ee
where
\bel{K_0}\sK^x_0[t,T]=\Big\{\xi(\cd)\equiv\sum_{k=1}^N\xi_k{\bf1}_{[\t_k,T]}(\cd)
\in\sK^x[t,T]\bigm|N\les{\bar\n(|x|)+1\over\ell_0},\q\sum_{k=1}^N|\xi_k|^\b\les {\bar\n(|x|)+1\over\a}\,\Big\}.\ee
Further, if {\rm(H1$'$)} is assumed, then the following holds:
\bel{V<*}0\les V(t,x)\les C(1+|x|^{\m+\d})+C_0e^{\b L(T-t)}\a(T)|x|^\b,\qq\forall(t,x)\in[0,T]\times\dbR^n,\ee
and in the definition of $\sK^x_0[t,T]$, $\bar\n(r)=C(1+r^{\m+\d})+C_0e^{\b L(T-t)}\a(T)r^\b$.

\ep

\it Proof. \rm For any $(t,x)\in[0,T]\times\dbR^n$, under the trivial impulse control $\xi_0(\cd)$, we have
$$|X(s;t,x,\xi_0(\cd))|\les e^{LT}(1+|x|),\qq t\les s< T.$$
Under (H1), making use of Lemma \ref{Lemma 2.6}, there is a $\xi_T\in K$ such that
$$X(T-0;t,x,\xi_0(\cd))+\xi_T\in\bar D,\qq|\xi_T|\les\n\big(|X(T-0;t,x,\xi_0(\cd))|\big),$$
Then we define
$$\h\xi(\cd)=\xi_T{\bf1}_{\{T\}}(\cd)\in\sK^x[t,T].$$
This is the impulse control that only makes one impulse at $T$ and make the state jump into $\bar D$. Clearly,
\bel{V(t,x)<}\ba{ll}
\ns\ds0\les V(t,x)\les J(t,x;\h\xi(\cd))=\int_t^Tg(s,X(s))ds+h\big(X(T-0)+\xi_T\big)
+\ell\big(T,X(T-0),\xi_T\big)\\
\ns\ds\les\int_t^TL\big(1+|X(s)|^{\m+\d}\big)ds+L(1+|X(T-0)+\xi_T|^{\m+\d})
+L+\a(T)|\xi_T|^\b\big)\les \bar\n(|x|),\ea\ee
for some continuous increasing function $\bar\n:[0,\i)\to[0,\i)$. This proves
\rf{V<}, which also leads to $\sD(V)=[0,T]\times\dbR^n$.

\ms

Next, let $(t,x)\in\sD(V)=[0,T]\times\dbR^n$, i.e., $V(t,x)<\infty$. Let $\xi(\cd)\in\sK^x[t,T]$ satisfy
$$V(t,x)+1\ges J(t,x;\xi(\cd))\ges N\ell_0+\sum_{k=1}^N\a|\xi_k|^\b.$$
Then
$$N\les{V(t,x)+1\over\ell_0}\les{\bar\n(|x|)+1\over\ell_0},$$
and
$$\ba{ll}
\ns\ds\sum_{k=1}^N|\xi_k|^\b\les{V(t,x)+1\over\a}\les {\bar\n(|x|)+1\over\a}.\ea$$
Hence, \rf{V=}--\rf{K_0} hold.

\ms

Now, let $\xi^\e(\cd)\in\sK^x_0[t,T]$ be a minimizing sequence for the cost functional $\xi(\cd)\mapsto J(t,x;\xi(\cd))$. Then, we may assume that
$$\lim_{\e\to0}\xi^\e(s)=\bar\xi(s)\equiv\sum_{k\ges1}\bar\xi_k{\bf1}_{[\bar\t_k,T]}
(\cd),\qq s\in[t,T],$$
with $t\les\bar\t_1<\bar\t_2<\cds$ and $\bar\xi_k\ne0$. This can be done as follows:
Let us begin with the sequence $(\t_1^\e,\xi_1^\e)$. We may assume that
$$\lim_{\e\to0}(\t_1^\e,\xi^\e_1)=(\t^0_1,\xi^0_1).$$
If $\xi^0_1\ne0$, we define
$$(\bar\t_1,\bar\xi_1)=(\t_1^0,\xi_1^0).$$
If $\xi_1^0=0$, we skip $(\t_1^0,\xi_1^0)$. By taking sub-subsequence, we may assume that
$$\lim_{\e\to0}(\t_2^\e,\xi_2^\e)=(\t_2^0,\xi_2^0).$$
If $\t_2^0>\t_1^0$ and $\xi_2^0\ne0$, we define
$$(\bar\t_2,\bar\xi_2)=(\t_2^0,\xi_2^0).$$
Otherwise, if $\xi_2^0=0$, we skip $(\t_2^0,\bar\xi_2^0)$; and if $\t_2^0=\t_1^0$, $\xi_2^0\ne0$, and $\xi_1^0+\xi_2^0=0$, we skip both $(\t_1^0,\xi_1^0)$ and $(\t_2^0,\xi_2^0)$; if $\t_2^0=\t_1^0$ and $\xi_1^0+\xi_2^0\ne0$, we redefine
$$(\bar\t_1,\bar\xi_1)=(\t_1^0,\xi_1^0+\xi_2^0).$$
Clearly, with such a procedure, we can complete the construction of $\bar\xi(\cd)$. By the convergence $\xi^\e(\cd)\to\bar\xi(\cd)$ pointwise, we see that actually the convergence is uniform. Then one also has the uniform convergence of $X^\e(\cd)\equiv X(\cd\,;t,x,
\xi^\e(\cd))$ to $\bar X(\cd)\equiv X(\cd\,;t,x,\bar\xi(\cd))$. It is ready to see that $(\bar X(\cd), \bar\xi(\cd))$ is an optimal pair.

\ms

Finally, if (H1$'$) holds, by Lemma \ref{Lemma 2.6}, we have from \rf{V(t,x)<}
that
$$\ba{ll}
\ns\ds0\les V(t,x)\les C(1+|x|^{\m+\d})+\a(T)C_0(1+|X(T-0)|^\b)\les C(1+|x|^{\m+\d})+C_0e^{\b L(T-t)}\a(T)|x|^\b,\ea$$
proving \rf{V<*}. \endpf

\ms

We see that the appearance of the term $C_0e^{\b L(T-t)}\a(T)|x|^\b$ is due to the terminal state constraint. It is possible that $\b>\m+\d$. In this case, the value function might grow with the same order as the impulse cost.

\subsection{Continuity of the value function}

In this section, we will establish the continuity of the value functions $V(\cd\,,\cd)$. Note that unlike the classical situation, when the terminal state constraint is presented, the value functions could be discontinuous. Also, some proper conditions will ensure the continuity of the value functions. To be convincing, let us first look at a simple example.

\bex{discontinuity} \rm Consider state equation (which is the same as that in Example \ref{Ex-2.4})
$$X(s)=x+(s-t)+\xi(s),\qq s\in[t,T].$$
The cost functional is defined by
$$J(t,x;\xi(\cd))=\sum_{k\ges1}\ell(\t_k,\wt X(\t_k-0),\xi_k),$$
with
$$\ell(t,x,\xi)=1+|\xi|.$$
Suppose $K=\dbR$, $D=(0,1)$. Let us consider Problem (C). For any $(t,x)\in[0,T]\times\dbR$ with $x+T-t\in[0,1]$,
we take $\xi(\cd)=\xi_0(\cd)$, the trivial impulse control. Hence,
$$V(t,x)=0,\qq\hb{ if }~x+T-t\in[0,1].$$
Now, if
$$x+T-t>1,$$
then, under $\xi_0(\cd)$, we have
$$X(T;t,x,\xi_0(\cd))=x+T-t>1.$$
Hence, during $[t,T]$ an impulse has to be made. The most economical impulse will be
$$\xi(\cd)=-(x+T-t-1){\bf1}_{[\t_1,T]}(\cd),$$
where the choice $\t_1\in[t,T]$ is irrelevant. Under such an impulse control, we have
$$J(t,x;\xi(\cd))=1+|x+T-t-1|=x+T-t.$$
Apparently, such an impulse control is optimal. Finally, if
$$x+T-t<0,$$
then we take
$$\xi(\cd)=-(x+T-t){\bf1}_{[\t_1,T]}(\cd),$$
with an arbitrary $\t_1\in[t,T]$. Again, this impulse control is optimal. With such a control, one has
$$J(t,x;\xi(\cd))=1+|x+T-t|=1-x-T+t.$$
Consequently,
\bel{V1}V(t,x)=\left\{\2n\ba{ll}
\ns\ds1-x-T+t,\qq\qq~~x+T-t<0,\\
\ns\ds0,\qq\qq\qq\qq\qq x+T-t\in[0,1],\\
\ns\ds x+T-t,\qq\qq\qq~ x+T-t>1.\ea\right.\ee
Clearly, this value function $V(\cd\,,\cd)$ is discontinuous (along the lines $x+T-t=0$ and $x+T-t=1$).

\ms

Now, we modify the cost functional as follows:
$$J(t,x;\xi(\cd))=h(X(T))+\sum_{k\ges1}\ell(\t_k,\wt X(\t_k-0),\xi_k),$$
with
$$h(x)=9\(x-{2\over5}\)^2,\qq x\in\dbR.$$
For any $X\in\dbR$ (a possible terminal state location), take $\xi\in K\equiv\dbR$ and look at the following:
$$h(X+\xi)+\ell(T,X,\xi)=9\(X+\xi-{2\over5}\)^2+1+|\xi|,$$
requiring $b\equiv X+\xi\in[0,1]$. This is the cost at the terminal time $T$ if the terminal state is $X$ and an impulse $\xi$ is made at $T$. Hence, let us consider the following function
$$F(b,X)=h(b)+\ell(T,X,b-X)=9\(b-{2\over5}\)^2+1+|b-X|,\qq b\in[0,1],$$
which will help us to decide whether we should make an impulse at $T$. For any given $X\in\dbR$, we want to find the minimum of $b\mapsto F(b,X)$ over $b\in[0,1]$. To this end, we first observe that
$$F_b(b,X)=\left\{\2n\ba{ll}
\ns\ds18\(b-{2\over5}\)+1=18b-{31\over5},\qq\hb{if }b>X,\\ [2mm]
\ns\ds18\(b-{2\over5}\)-1=18b-{41\over5},\qq\hb{if }b<X,\ea\right.$$
Hence,
$$\left\{\2n\ba{ll}
\ns\ds F_b\({31\over90},X\)=0,\qq X<{31\over90}\equiv b_0,\\ [3mm]
\ns\ds F_b\({41\over90},X\)=0,\qq X>{41\over90}\equiv b_1.\ea\right.$$
Clearly,
$$0<b_0\equiv{31\over90}<{2\over5}<b_1\equiv{41\over90}<1.$$
Further, for $X\in(b_0,b_1)$, we have
$$\left\{\2n\ba{ll}
\ns\ds F_b(b,X)<0,\qq b\in(b_0,X),\\
\ns\ds F_b(b,X)>0,\qq b\in(X,b_1).\ea\right.$$
Hence, for $X\in(b_0,b_1)$,
$$\min_{b\in[0,1]}F(b,X)=F(X,X)=9\(X-{2\over5}\)^2+1.$$
To summarize, we have
$$\min_{b\in[0,1]}F(b,X)=\left\{\2n\ba{ll}
\ns\ds F(b_0,X)=9\({31\over90}-{2\over5}\)^2+1+{31\over90}-X
={247\over180}-X,\qq X<b_0,\\ [2mm]
\ns\ds F(X,X)=9\(X-{2\over5}\)^2+1,\qq X\in[b_0,b_1],\\ [2mm]
\ns\ds F(b_1,X)=9\({41\over90}-{2\over5}\)^2+1+X-{41\over90}
={103\over180}+X,\qq X>b_1,\ea\right.$$
and $\ds X\mapsto\min_{b\in[0,1]}F(b,X)$ is continuous. Note that
\bel{h(0)>}h(0)={36\over25}>{247\over180}=\min_{b\in[0,1]}F(b,0)=\min_{\xi\in[0,1]}
\[h(\xi)+\ell(T,0,\xi)\],\ee
and
\bel{h(1)>}h(1)={81\over25}>{283\over180}=\min_{b\in[0,1]}F(b,1)=\min_{\xi\in[0,1]}
\[h(1+\xi)+\ell(T,1,\xi)\],\ee
Now, we look at the equation
$$h(X)=\min_{b\in[0,1]}F(b,X),$$
which give the point $X$ at which there is no difference if the best impulse is made or no impulse is made. A direct check shows that the above does not have solutions in $[b_0,b_1]$. Now, on $(0,b_0)$, we solve
$$9\(X-{2\over5}\)^2={247\over180}-X,$$
whose unique solution is
$$X={1\over90}\equiv a_0.$$
On $(b_1,1)$, we solve
$$9\(X-{2\over5}\)^2={103\over180}+X,$$
whose unique solution is
$$X={71\over90}\equiv a_1.$$
The above tells us that (recalling $b=X+\xi$)
$$\left\{\2n\ba{ll}
\ns\ds h(X)<\min_{X+\xi\in[0,1]}\[h(X+\xi)+\ell(T,X,\xi)\],\qq X\in(a_0,a_1),\\
\ns\ds h(X)\ges\min_{X+\xi\in[0,1]}\[h(X+\xi)+\ell(T,X,\xi)\],\qq X\in\dbR\setminus
(a_0,a_1).\ea\right.$$
We see that
$$0<a_0\equiv{1\over90}<b_0\equiv{31\over90}<b_1\equiv{41\over90}<a_1\equiv
{71\over90}<1.$$
The above means that if the terminal state $X(T-0)\in(a_0,a_1)$, we should not make an impulse at $T$, and if $X(T-0)\in\dbR\setminus(a_0,a_1)$, we should make an impulse as follows:
\bel{xi_1}\xi_1=\left\{\1n\ba{ll}
\ds{31\over90}-X(T-0),\qq X(T-0)<a_0={1\over90},\hb{ or }x+T-t<a_0,\\ [3mm]
\ds{41\over90}-X(T-0),\qq X(T-0)>a_1={71\over90},\hb{ or }x+T-t>a_1,\ea\right.\ee
so that either $X(T)=b_0\equiv{31\over90}$ or $X(T)=b_1\equiv{41\over90}$.
Combining the above analysis, we obtain the value function
\bel{V2}V(t,x)=\left\{\2n\ba{ll}
\ns\ds9\(x+T-t-{2\over5}\)^2,\qq x+T-t\in\[{1\over90},{71\over90}\]\equiv[a_0,a_1],\\ [2mm]
\ns\ds{247\over180}-(x+T-t),\qq~x+T-t<{1\over90}\equiv a_0,\\ [2mm]
\ns\ds{103\over180}+x+T-t,\qq\q x+T-t>{71\over90}\equiv a_1,\ea\right.\ee
which is continuous.

\ms

Now, let $K=[0,\infty)$ and $D=(0,1)$. Then from Example \ref{Ex-2.4}, we see that
$$\sD(V)=\big\{(t,x)\in[0,T]\times\dbR\bigm|x+T-t\les1\big\},$$
and only positive impulses can be made. Hence, by looking above computation, we see that
if  $X(T-0;t,x,\xi_0(\cd))<a_0$, we could make an impulse; for all other cases, we could not/should not make impulses. Therefore,
\bel{V3}V(t,x)=\left\{\2n\ba{ll}
\ds{247\over180}-(x+T-t),\qq~x+T-t<{1\over90},\\ [2mm]
\ns\ds9\(x+T-t-{2\over5}\)^2,\qq x+T-t\in\[{1\over90},1\],\\ [2mm]
\ns\ds+\infty,\qq\qq\qq\qq~x+T-t>1.\ea\right.\ee
This value function is continuous over $\sD(V)=\cl{\sD(V)}$ which is a closed set.

\ms

Further, let $K=[0,\infty)$ and $D=(0,\infty)$. The feature is that the state $X(T-0)$ will be either in $\bar D$, or, it can always be pulled back to $\bar D$ by an admissible impulse. Therefore, $\sD(V)=[0,T]\times\dbR$. When $X(T-0)<a_0={1\over90}$, an impulse is necessary to meet the terminal state constraint, or to make the total cost smaller. Hence,
\bel{V4}V(t,x)=\left\{\2n\ba{ll}
\ns\ds{247\over180}-(x+T-t),\qq~x+T-t<{1\over90},\\ [2mm]
\ns\ds9\(x+T-t-{2\over5}\)^2,\qq x+T-t\ges{1\over90}.\ea\right.\ee
This function is continuous as well.

\ex

The above example shows that when the terminal cost function $h(\cd)$ and the impulse cost are compatible, one could get the continuity of the value function $V(\cd\,,\cd)$. In the above example, A careful observation shows that when the terminal state gets close to the boundary $\pa D$ of the constraint set $D$ from inside, an impulse should be made to reduce the cost. This essentially eliminates the possible jumps of the best costs between the terminal state $X(T-0)$ being close to the boundary $\pa D$ from outside and from inside of $D$. On the other hand, due to the terminal constraint, the value $h(x)$ of $h(\cd)$ for $x\in\dbR^n\setminus\bar D$ is irrelevant to our problem. We now would like to present general results.

\bt{continuity} \sl Let {\rm(H1)--(H4)} hold and \rf{D-K=R} be satisfied. Suppose
\bel{inf<}\inf_{\xi\in K,x+\xi\in D}\[h(x+\xi)+\ell(T,x,\xi)\]<h(x),\qq\forall x\in\dbR^n\setminus D.\ee
Then there exists a continuous increasing function $\h\n:[0,\i)\to[0,\i)$ such that
\bel{V-V<L}\left\{\2n\ba{ll}
\ds|V(t,x)-V(t,\h x)|\les C(1+|x|^\m\vee|\h x|^\m)|x-\h x|^\d,\qq\forall x,\h x\in\dbR^n,~|x-\h x|\hb{ small},\\
\ns\ds|V(t,x)-V(\,\h t,x)|\les\h\n\big(|x|\vee|\h x|\big)|t-\h t\,|\qq\forall t,\h t\in[0,T].\ea\right.\ee
In the case that {\rm(H1$'$)} holds,
\bel{V-V<L*}|V(t,x)-V(\,\h t,\h x)|\les C\big(1+|x|^\m\vee|\h x|^\m\big)(|t-\h t\,|
+|x-\h x|^\d),\q\forall(t,x),(\,\h t,\h x)\in[0,T]\times\dbR^n~|x-\h x|\hb{ small}.\ee

\et

\it Proof. \rm Let $(t,x)\in[0,T]\times\dbR^n$. From Proposition \ref{Prop-3.2}, there exists an optimal impulse control $\bar\xi(\cd)\in\sK^x_{0}[t,T]$. Due to condition \rf{inf<}, we claim that
$\bar X(T)\in D$. In fact, if $\bar X(T)\in\pa D$, then there exists a $\z\in K$ such that
$$\bar X(T)+\z\in D,\qq h(\bar X(T)+\z)+\ell(T,\bar X(T),\z)<h(\bar X(T)).$$
Thus, by letting
$$\h\z(\cd)=\bar\xi(\cd)+\z{\bf1}_{\{T\}},$$
we have
$$J(t,x;\h\z(\cd))=J(t,x;\bar\xi(\cd))+h(\bar X(T)+\z)+\ell(T,\bar X(T),\z)-h(\bar X(T))<J(t,x;\bar\xi(\cd))=V(t,x),$$
contradicting the optimality of $\bar\xi(\cd)$. Hence, we may assume that $\bar X(T)\in D$. Now, for any $\h x\in\dbR^n$, let $\h X(\cd)=X(\cd\,;t,\h x,\bar\xi(\cd))$, we have
$$|\bar X(s)-\h X(s)|\les e^{L(T-t)}|x-\h x|,\qq t\les s\les T.$$
Recalling that $D$ is open, for $|x-\h x|$ small, one sees that $\bar\xi(\cd)\in\sK^{\h x}[t,T]$. Consequently, making use of Propositions \ref{state X} and \ref{Prop-3.2}, together with the Lipschitz continuity of $x\mapsto\ell(t,x,\xi)$, we have (noting $\bar\xi(\cd)\in\sK^x_0[t,T]$)
$$\ba{ll}
\ns\ds V(t,\h x)\les J(t,\h x;\bar\xi(\cd))\les J(t,x;\bar\xi(\cd))+L\int_t^T
\(1+|\bar X(r)|^\m\vee|\h X(r)|^\m\)|\bar X(r)-\h X(r)|^\d dr\\
\ns\ds\qq\qq\qq\qq\qq\qq+L\(1+|\bar X(T)|^\m\vee|\h X(T)|^\m\)|\bar X(T)-\h X(T)|^\d\\
\ns\ds\qq\qq\qq\qq\qq\qq+\sum_{k\ges1}|\ell(\bar\t_k,\bar X(\bar\t_k-0),\bar\xi_k)-\ell(\bar\t_k,\h X(\bar\t_k-0),\bar\xi_k)|\\
\ns\ds\qq\q\les V(t,x)+C(1+|x|^\m\vee|\h x|^\m)|x-\h x|^\d.\ea$$
for some continuous increasing function $\h\n:[0,\i)\to(0,\i)$.
By symmetry, we obtain \rf{V-V<L}.

\ms

Next, let $0\les t<\h t\les T$. Let $\bar\xi(\cd)\in\sK_0^x[t,T]$ be optimal for the initial pair $(t,x)$. Let
$$t\les\bar\t_1<\bar\t_2<\cds<\bar\t_{k_0}\les\h t<\bar\t_{k_0+1}.$$
Define
$$\h\xi(\cd)=\sum_{k=1}^{k_0}\bar\xi_k{\bf1}_{[\hat t,T]}(\cd)+\sum_{k\ges k_0+1}\bar\xi_k{\bf1}_{[\bar\t_k,T]}(\cd).$$
Denote $\h X(\cd)=X(\cd\,;\h t,x,\h\xi(\cd))$. Then
$$\ba{ll}
\ns\ds J(\,\h t,x;\h\xi(\cd))=\int_{\h t}^Tg(r,\h X(r))dr+h(\h X(T))
+\ell\(\h t,x;\sum_{\bar\t_k\les\h t}\bar\xi_k\)+\sum_{\t_k>\h t}\ell(\bar\t_k,\h X(\bar\t_k-0),\bar\xi_k)\\
\ns\ds\les J(t,x;\bar\xi(\cd))+\int_t^{\h t}|g(r,X(r))|dr+\int_{\h t}^T|g(r,\h X(r))
-g(r,X(r))|dr+|h(\h X(T))-h(X(T))|\\
\ns\ds\qq+\[\ell\(\h t,x;\sum_{k=1}^{k_0}\bar\xi_k\)-\sum_{k=1}^{k_0}\ell(\bar\t_k,X(\bar\t_k-0),\bar
\xi_k)\]+\sum_{k>k_0}\[\ell(\bar\t_k,\h X(\t_k-0),\bar\xi_k)-\ell(\bar\t_k,X(\bar\t_k-0),\bar\xi_k)\].\ea$$
Note that
$$\ell\(\h t,x;\sum_{k=1}^{k_0}\bar\xi_k\)\les\sum_{k=1}^{k_0}\ell\(\h t,x+\sum_{i=1}^{k-1}
\bar\xi_i,\bar\xi_k\)\les\sum_{k=1}^{k_0}\ell\(\bar\t_k,x+\sum_{i=1}^{k-1}\bar\xi_i,
\bar\xi_k\).$$
On the other hand,
$$\ba{ll}
\ns\ds|X(\t_1-0)-x|\les\int_t^{\t_1}|f(r,X(r))|dr\les L\int_t^{\t_1}(1+|X(r)|)dr\\
\ns\ds\les L\int_t^{\t_1}\(1+e^{L(r-t)}\big(1+|x|\big)\)dr\les C(1+|x|)(\t_1-t).\ea$$
Next,
$$\ba{ll}
\ns\ds|X(\bar\t_2-0)-x-\xi_1|\les\int_{\bar\t_1}^{\bar\t_2}|f(r,X(r))|dr\les L\int_{\bar\t_1}^{\bar\t_2}\big(1+|X(r)|\big)dr\\
\ns\ds\les L\int_{\bar\t_1}^{\bar\t_2}\(1+e^{L(r-t)}(1+|x|)+e^{L(r-\bar\t_1)}|\bar\xi_1|\)dr\les
C(1+|x|+|\bar\xi_1|)(\bar\t_2-\bar\t_1).\ea$$
By induction, we see that
$$\Big|X(\bar\t_k-0)-x-\sum_{i=1}^{k-1}\bar\xi_i\Big|\les C\(1+|x|+\sum_{i=1}^{k-1}|\bar\xi_i|\)(\bar\t_k-t),\qq1\les k\les k_0.$$
Also, for any $s>\hat t$,
$$\ba{ll}
\ns\ds|\h X(s)-X(s)|\les e^{L(s-\h t\,)}\Big|x+\sum_{i=1}^{k_0}\bar\xi_k-X(\hat t\,)\Big|\les C\(1+|x|+\sum_{i=1}^{k_0}|\bar\xi_i|\)(\h t-t).\ea$$
Consequently, noting that $\bar\xi(\cd)\in\sK^x_0[t,T]$,
$$\ba{ll}
\ns\ds J(\,\h t,x;\h\xi(\cd))\les J(t,x;\bar\xi(\cd))+\int_t^{\h t}L(1+|X(r)|^\d)dr+L\int_{\h t}^T|X(r)-\h X(r)|dr+L|X(T)-\h X(T)|\\
\ns\ds\qq\qq\qq\qq\qq+L\sum_{k=1}^{k_0}\Big|x+\sum_{i=1}^{k-1}\bar\xi_i
-X(\bar\t_k-0)\Big|+L\sum_{k>k_0}|\h X(\bar\t_k-0)-X(\bar\t_k-0)|\\
\ns\ds\qq\qq\qq\les V(t,x)+C\(1+|x|+\sum_{k\ges1}|\bar\xi_k|\)(\,\h t-t)\les V(t,x)+\wt\n(|x|)(\,\h t-t),\ea$$
for some $\wt\n(|x|)$. Finally, let $\h t<t$. Then we extend $\bar\xi(\cd)$ on $[t,T]$ to $\h\xi(\cd)$ on $[\,\h t,T]$ trivially. One has
$$\ba{ll}
\ns\ds V(\,\h t,x)\les J(\,\h t,x;\h\xi(\cd))\les J(t,x;\bar\xi(\cd))
+\int_{\h t}^t|g(r,\h X(r))|dr+\int_t^T|g(r,\h X(r))-g(r,X(r))|dr\\
\ns\ds\qq\qq+|h(\h X(T))-h(X(r))|+\sum_{k\ges1}|\ell(\bar\t_k,\h X(\bar\t_k-0),
\xi_k^\e)-\ell(\bar\t_k,X(\bar\t_k-0),\bar\xi_k^\e)|\\
\ns\ds\qq\les V(t,x)+\int_{\h t}^tL(1+|\h X(r)|)dr+L\int_t^T|X(r)-\h X(r)|dr\\
\ns\ds\qq\qq+L|X(T)-\h X(T)|+\sum_{k\ges1}
|X(\bar\t_k-0)-\h X(\bar\t_k-0)|\\
\ns\ds\qq\les V(t,x)+C\(1+|x|+\sum_{k\ges1}|\bar\xi_k|\)
(t-\h t\,)\les\wt\n(|x|)(t-\h t\,).\ea$$
This completes the proof of \rf{V-V<L}. Finally, in the case that (H1$'$) holds, our conclusion follows from the above arguments, together with Proposition \ref{Prop-3.2}. \endpf

\ms

We see that due to the appearance of the terminal constraint, the H\"older continuity of the map $x\mapsto V(t,x)$ is only locally.

\section{Dynamic Programming Principle and HJB Quasi-Variational Inequality}

In this section, we are going to establish Bellman's principle of optimality for our Problem (C). Then the corresponding HJB equation, which is a quasi-variational inequality, for
the value function $V(\cd\,,\cd)$ will be derived. For convenience, in what follows, we will keep assumptions (H1$'$), (H2)--(H4) and \rf{inf<}.

\bt{DPP} \sl Let $V(\cd\,,\cd)$ be the value function of Problem (C). Then for any $(t,x)\in[0,T)\times\dbR^n$, the following principle of optimality holds:
\bel{optimality-1}
V(t,x)\les\min_{\xi\in K}\big\{V(t,x+\xi)+\ell(t,x,\xi)\big\}\equiv\BN[V](t,x),\qq\forall(t,x)\in[0,T)\times \dbR^n,\ee
\bel{optimality-2}
V(t,x)\les\int_t^{\h t}g(r,X(r;t,x,\xi_0(\cd)))dr+V(\,\h t,X(\h t;t,x,\xi_0(\cd))),\q \forall0\les t\les\h t\les T,~x\in\dbR^n.\ee
Furthermore, if the strict inequality holds in \rf{optimality-1}, then there exists
a $\bar t\in(t,T]$ such that
$$V(t,x)=\int_t^{\h t}g(r,X(r;t,x,\xi_0(\cd)))dr+V(\,\h t,X(\h t;t,x,\xi_0(\cd))),\qq  0\les t\les\h t<\bar t\les T,~x\in\dbR^n.$$
For $x\in\dbR^n$, it holds that
\bel{V(T)}V(T,x)=h^D(x)\equiv\left\{\2n\ba{ll}
\ds\min\big\{h(x),\BN^D[h](x)\big\},\qq x\in\bar D,\\ [2mm]
\ns\ds\BN^D[h](x),\qq\qq\qq\q x\in\dbR^n\setminus\bar D,\ea\right.\ee
where
\bel{h BN}\BN^D[h](x)=\min_{\xi\in K,x+\xi\in\bar D}\big\{h(x+\xi)+\ell(T,x,\xi)\big\},\qq x\in\dbR^n.\ee
\et

\it Proof. \rm First of all, for $x\in\dbR^n$, we clearly have \rf{V(T)}. Next, let $(t,x)\in[0,T)\times\dbR^n$ and $\z\in K$. For any $\xi(\cd)\in\sK^{x+\z}[t,T]$, we see that
$$\h\xi(\cd)=\z{\bf1}_{[t,T]}(\cd)+\xi(\cd)\in\sK^x[t,T].$$
Thus,
$$V(t,x)\les J(t,x;\h\xi(\cd))=\ell(t,x,\z)+J(t,x+\z;\xi(\cd)).$$
Consequently,
$$V(t,x)\les \ell(t,x,\z)+V(t,x+\z),\qq\forall\z\in K.$$
Therefore,
$$V(t,x)\les\min_{\z\in K}\{V(t,x+\z)+\ell(t,x,\z)\}=\BN[V](t,x).$$
On the other hand, for any $0\les t<\h t\les T$, we take any $\xi(\cd)\in\sK^{X(\h t;t,x,\xi_0(\cd))}[\,\h t,T]$. Extend it to $\h\xi(\cd)\in\sK^x[t,T]$ in the way that
no impulses are made on $[t,\h t\,)$. Then
$$V(t,x)\les J(t,x;\h\xi(\cd))=\int_t^{\h t}g(r,X_{t,x}(r))dr+J(\,\h t,X(\,\h t;t,x,\xi_0(\cd));\xi(\cd)).$$
Consequently,
$$V(t,x)\les\int_t^{\h t}g(r,X(r;t,x,\xi_0(\cd)))dr+V(\,\h t,X(\,\h t;t,x,\xi_0(\cd))).$$

Finally, we assume that
\bel{assu}V(t,x)<\BN[V](t,x)\equiv\min_{\z\in K}\big\{V(t,x+\z)+\ell(t,x,\z)\big\}.\ee
For any $\e>0$, there exists an impulse control $\ds\xi^\e(\cd)=\sum_{i\ges1}\xi_i^\e{\bf1}_{[\t_i^\e,T]}(s)\in\sK^x[t,T]$ such that
$$J(t,x;\xi^\e(\cd))\les V(t,x)+\e.$$
If $\t_1^\e=t$, then
\bel{}V(t,x)+\e\ges J(t,x+\xi_1^\e;\h\xi^\e(\cd))+\ell(t,x,\xi_1^\e)
\ges V(t,x+\xi_1^\e)+\ell(t,x,\xi_1^\e),\ee
where
$$\h\xi^\e(\cd)=\sum_{i\ges1}\h\xi_i^\e{\bf1}_{[\hat\t_i^\e,T]}(s);\qq\h\xi_i^\e
=\xi_{i+1}^\e,\q\hat\t_i^\e=\t_{i+1}^\e,\q i\ges1.$$
This is contradicting \rf{assu}. Hence, $\t_1^\e>t$. We further claim that there exists a $\bar t>t$ such that $\t_1^\e\ges\bar t$ for all $\e>0$ small. If this is not the case, then for some $\e\downarrow0$, $\t_i^\e\downarrow t$. Thus,
\bel{}\ba{ll}
\ns\ds V(t,x)+\e\ges J(t,x;\xi^\e(\cd))\\
\ns\ds=\int_t^{\t_1^\e}\3n g(r,X(s;t,x,\xi_0(\cd))ds+\ell(\t_1^\e,X(\t^\e_1-0;t,x,\xi_0(\cd)),
\xi_1^\e)+J(\t_1^\e,
X(\t_1^\e-0;t,x,\xi_0(\cd))+\xi_1^\e;\h\xi^\e(\cd))\\
\ns\ds\ges\int_t^{\t_1^\e}g(r,X(s;t,x,\xi_0(\cd)))ds+\ell(\t_1^\e,X(\t_1^\e-0;t,x\xi_0(\cd)),
\xi_1^\e)+V\big(\t_1^\e,X(\t_1^\e-0;t,x,\xi_0(\cd))\big).\ea\ee
Letting $\e\downarrow 0$ and (we may assume that) $\xi^\e_1\rightarrow\xi_1$, we have
$$V(t,x)\ges V(t,x+\xi_1)+\ell(t,x,\xi_1).$$
which is a contradiction again. Therefore, we get the existence of $\bar t>t$. Then for any $\h t<\bar t$,
\bel{}\ba{ll}
\ns\ds V(t,x)+\e\ges\int_t^{\h t}g(r,X(r;t,x,\xi_0(\cd)))dr+J(\,\h t,X(\,\h t;t,x,\h\xi^\e)\\
\ns\ds\qq\qq\q\ges\int_t^{\,\h t}g(r,X(r;t,x,\xi_0(\cd)))ds+V(\,\h t,X(\,\h t;t,x,\xi_0(\cd))).\ea\ee
Combing this with \rf{optimality-2}, the proof is complete.
\endpf

\ms

The following proposition is about the properties of $\BN[\,\cd\,]$ and $\BN^D[\,\cd\,]$.

\bp{BN} \sl The maps $(t,x)\mapsto\BN[V](t,x)$ and $x\mapsto\BN^D[h](x)$ are continuous. Moreover, if $(t_0,x_0)\in[0,T)\times\dbR^n$ such that
\bel{V=BN}V(t_0,x_0)=\BN[V](t_0,x_0)=V(t_0,x_0+\xi_0)+\ell(t_0,x_0,\xi_0),\ee
for some $\xi_0\in K$, then
\bel{V<BN}\BN[V](t_0,x_0+\xi_0)-V(t_0,x_0+\xi_0)\ges\d_0.\ee
Similarly, if $x_0\in\dbR^n$ such that
\bel{}h(x_0)=\BN^D[h](x_0)=h(x_0+\xi_0)+\ell(T,x_0,\xi_0),\ee
for some $\xi_0\in K$, then
\bel{}\BN^D[h](x_0+\xi_0)-h(x_0+\xi_0)\ges\d_0.\ee

\ep

\it Proof. \rm First of all, suggested by \rf{V=}--\rf{K_0}, for any $x\in\dbR^n$, we may introduce
\bel{K_0^x}K_0^{|x|}=\Big\{\xi\in K\bigm||\xi|\les\({\bar\n(|x|)+1\over\a}\)^{1\over\b}\Big\}.\ee
Then \rf{V=}--\rf{K_0} tells us that for given initial pair $(t,x)\in[0,T]\times\dbR^n$, we may restrict ourselves to the impulse controls with the impulse vector taken from
$K_0^{|x|}$. Now, for any $(t,x),(\,\h t,\h x)\in[0,T)\times\dbR^n$, and $\xi\in
K_0^{|x|\vee|\h x|}$,
$$\big|V(t,x+\xi)+\ell(t,x,\xi)-[V(\,\h t,\h x+\xi)+\ell(\,\h t,\h x,\xi)]\big|\les\big[\h\n(|x|\vee|\h x|)+L\big]\big(|t-\h t\,|+|x-\h x|\big).$$
Note that in the definition of $\BN[\,\cd\,]$, $\xi\in K$ is uniform in $(t,x)$. Thus, we have
\begin{equation}\label{inf-inf-1}
|\BN[V](t,x)-\BN[V](\,\h t,\h x)|\les\big[\h\n(|x|\vee|\h x|)+L\big]\big(|t-\h t\,|+|x-\h x|\big).
\end{equation}
Next, if for some $(t_0,x_0)\in[0,T)\times\dbR^n$, \rf{V=BN} holds for some $\xi_0\in
K_0^{|x_0|}$, then we claim that \rf{V<BN} holds. In fact, for any $\xi_1\in K_0^{|x_0+\xi_0|}$
$$\ba{ll}
\ns\ds V(t_0,x_0+\xi_0+\xi_1)+\ell(t_0,x_0+\xi_0,\xi_1)-V(t_0,x_0+\xi_0)\\
\ns\ds=V(t_0,x_0+\xi_0+\xi_1)+\ell(t_0,x_0+\xi_0,\xi_1)+\ell(t_0,x_0,\xi_0)
-\[V(t_0,x_0+\xi_0)+\ell(t_0,x_0,\xi_0)\]\\
\ns\ds=\1n V(t_0,x_0\1n+\1n\xi_0\1n+\1n\xi_1)\1n+\1n\ell(t_0,x_0,\xi_0\1n+\1n\xi_1)\1n-\1n
V(t_0,x_0)\1n+\1n\[\ell(t_0,x_0,\xi_0)\1n
+\1n\ell(t_0,x_0\1n+\1n\xi_0,\xi_1)\1n-\1n\ell(t_0,x_0,\xi_0\1n+\1n\xi_1)\]\\
\ns\ds\ges\ell(t_0,x_0,\xi_0)
+\ell(t_0,x_0+\xi_0,\xi_1)-\ell(t_0,x_0,\xi_0+\xi_1)\ges\d_0.\ea$$
 This proves \rf{V<BN}.

\ms

Now, we look at $\BN^D[\,\cd\,]$. By \rf{K_0^x}, we can redefine (compare with \rf{h BN})
$$\BN^D[h](x)=\min_{\xi\in K_0^x,x+\xi\in D}\big\{h(x+\xi)+\ell(T,x,\xi)\big\},\qq x\in\dbR^n.$$
We now show that $x\mapsto\BN^D[h](x)$ is continuous. For any given $x\in\dbR^n$, there exists a $\xi\in K_0^x\subset K$ such that
$$\BN^D[h](x)=h(x+\xi)+\ell(T,x,\xi),\qq x+\xi\in\bar D.$$
We claim that $x+\xi\in D$ (not on the boundary $\pa D$ of $\bar D$). In fact, if $x+\xi\in\pa D$, then by \rf{inf<}, there exists a $\xi_1\in K$ such that
$$\ba{ll}
\ns\ds\BN^D[h](x)=h(x+\xi)+\ell(T,x,\xi)>h(x+\xi+\xi_1)+\ell(T,x+\xi,\xi_1)
+\ell(T,x,\xi)\\
\ns\ds\qq\qq>h(x+\xi+\xi_1)+\ell(T,x,\xi+\xi_1)\ges\BN^D[h](x),\ea$$
a contradiction. Hence, by the openness of $D$, there exists a $\d>0$ such that $\sO_\d(x+\xi)\subseteq D$. Then for $y\in\sO_\d(x)$, $\xi\in D-\{y\}$ and thus
$$\ba{ll}
\ns\ds\BN^D[h](x)=h(x+\xi)+\ell(T,x,\xi)\ges h(y+\xi)+\ell(T,y,\xi)-\n(|x|\vee|y|,|x-y|)\\
\ns\ds\qq\qq\q\ges\BN^D[h](y)-\n(|x|\vee|y|,|x-y|),\ea$$
for some continuous function $\n:[0,\i)\times[0,\i)\to[0,\i)$ with $\n(r,0)=0$, for any $r\ges0$. Switch the positions of $x$ and $y$, we obtain the continuity of $x\mapsto\BN^D[h](x)$. The proof of last conclusion is similar to the case of $\BN[V](\cd\,,\cd)$, with the restriction that $x_0+\xi_0,x_0+\xi_0+\xi_1\in\bar D$ and
$t_0=T$. \endpf

\ms

The above result leads to the following Hamilton-Jacobi-Bellman equation for the value function $V(\cd\,,\cd)$. The proof is standard.

\bt{HJB} \sl Suppose that the value function $V(\cd\,,\cd)\in C^1([0,T]\times\dbR^n)$. Then $V(\cd\,,\cd)$ satisfies the following {\it HJB quasi-variational inequality}:
\bel{QI}\left\{\2n\ba{ll}
\ns\ds V_t(t,x)+\lan V_x(t,x),f(t,x)\ran+g(t,x)\ges0,\qq\BN[V](t,x)-V(t,x)\ges0,\q(t,x)\in[0,T]\times\dbR^n,\\
\ns\ds\(V_t(t,x)+\lan V_x(t,x),f(t,x)\ran+g(t,x)\)\(\BN[V](t,x)-V(t,x)\)=0,\q(t,x)\in[0,T]\times\dbR^n,\\
\ns\ds V(T,x)=h^D(x),\qq x\in\dbR^n,\ea\right.\ee
which can also be written as
\bel{HJB-1}\left\{\ba{ll}
\ds\min\Big\{V_t(t,x)+\lan V_x(t,x),f(t,x)\ran+g(t,x),\BN[V](t,x)-V(t,x)\Big\}=0,\q
(t,x)\in[0,T)\times\dbR^n,\\
\ns\ds V(T,x)=h^D(x),\qq x\in\dbR^n.\ea\right.\ee
\et

From the previous sections, we see that under (H1$'$), (H2)--(H4) and \rf{inf<}, the value function $V(\cd\,,\cd)\in C([0,T]\times\dbR^n)$. However, it is known that the value function might not be $C^1([0,T]\times\dbR^n)$ in general. Therefore, the above is a formal result. For \rf{HJB-1}, inspired by the viscosity solution notion introduced by Crandall--Lions \cite{Crandall-Lions 1983}, Barles introduced the following corresponding notion \cite{Barles 1985, Barles 1985b}, which has been modified here for our impulse control problem in finite time horizon (see \cite{Tang-Yong 1993}).

\bde{Viscosity solution} \rm A continuous function $V(\cd\,,\cd)$ is called a {\it viscosity sub-solution} of HJB quasi-variational inequality \rf{HJB-1} if
\bel{VT)<}V(T,x)\les h^D(x),\qq x\in\dbR^n,\ee
and for any function $\f\in C^1([0,T]\times\dbR^n)$ such that $V(\cd\,,\cd)-\f(\cd\,,\cd)$ attains a local maximum at $(t,x)\in[0,T)\times\dbR^n$, it holds
\bel{min>0}\min\Big\{\f_t(t,x)+\lan\f_x(t,x),f(t,x)\ran+g(t,x),\BN[V](t,x)-V(t,x)\Big\}\ges0.\ee
A continuous function $V(\cd\,,\cd)$ is called a {\it viscosity super-solution} of HJB quasi-variational inequality \rf{HJB-1} if
\bel{V(T)>}V(T,x)\ges h^D(x),\qq x\in\dbR^n,\ee
and for any function $\f\in C^1([0,T]\times\dbR^n)$ such that $V(\cd\,,\cd)-\f(\cd\,,\cd)$ attains a local minimum at $(t,x)\in[0,T)\times\dbR^n$, it holds
\bel{min<}\min\Big\{\f_t(t,x)+\langle\f_x(t,x),f(t,x)\rangle+g(t,x),\BN[V](t,x)-V(t,x)\Big\}\les0.\ee
A continuous function $V(\cd\,,\cd)$ is called a {\it viscosity solution} of the HJB quasi-variational inequality \rf{HJB-1} if it is both viscosity super-solution and viscosity sub-solution.
\ede


%
%
%


We now state the following result whose proof is (almost) standard (see \cite{Barles 1985, Barles 1985b, Tang-Yong 1993}).

\bt{viscosity-1} \sl Let {\rm(H1$'$), (H2)--(H4)} and \rf{inf<} hold, then the value function $V(\cd\,,\cd)$ is
a viscosity solution to the HJB quasi-variational inequality \rf{HJB-1}.

%
%
%

\et

We have seen that Problem (C) admits optimal impulse control (see Proposition \ref{Prop-3.2}). It is almost standard that, via the value function, an optimal impulse control can also be constructed. We omit the details here.

\section{An Optimal Impulse Control Problem without Terminal State Constraint}

HJB QVI \rf{HJB-1} suggests us to introduce the following modified cost functional:
\bel{wtJ}\wt J(t,x;\xi(\cd))=\int_t^Tg(s,X(s))ds+h^D(X(T))+\sum_{k\ges1}\ell(\t_k,\wt X(\t_k-0),\xi_k),\ee
and consider the problem with the above cost functional without terminal state constraint, call it Problem ($\wt{\rm C}$). Let the value function be $\wt V(\cd\,,\cd)$. The question is whether $\wt V(\cd\,,\cd)=V(\cd\,,\cd)$? If this is true, then since $\wt V(\cd\,,\cd)$ can be characterized as the unique viscosity solution to \rf{HJB-1}, we indirectly obtain a characterization of the value function $V(\cd\,,\cd)$. In this section, we discuss this issue.

\ms

We recall that
\bel{N^D}\BN^D[h](x)=\min_{\xi\in K,x+\xi\in\bar D}\(h(x+\xi)+\ell(T,x,\xi)\),\qq x\in\dbR^n,\ee
and (for convenience, we denote the value function by $V^D(t,x)$)
\bel{V(T)=..}V^D(T,x)=h^D(x)\equiv\left\{\2n\ba{ll}
\ds\min\Big\{h(x),\BN^D[h](x)\Big\},\qq x\in\bar D,\\
\ns\ds\BN^D[h](x),\qq\qq\qq\qq x\in\dbR^n\setminus\bar D.\ea\right.\ee
Taking $D=\dbR^n$, the above is reduced to the case of no terminal state constraint. In other words, if there is no terminal state constraint, then the terminal value of the value function is given by the following:
\bel{V(T)=...}V^{\dbR^n}(T,x)=h^{\dbR^n}(x)\equiv\min\Big\{h(x),\BN^{\dbR^n}[h](x)\Big\},\qq x\in\dbR^n,\ee
with
\bel{N^R}\BN^{\dbR^n}[h](x)=\min_{\xi\in K}\(h(x+\xi)+\ell(T,x,\xi)\),\qq x\in\dbR^n.\ee
According to the above, we have
$$\wt V(T,x)=(h^D)^{\dbR^n}(x)\equiv\min\Big\{h^D(x),\BN^{\dbR^n}[h^D](x)\Big\},\qq
x\in\dbR^n.$$
Let us calculate the following: (remember $D\ne\dbR^n$ and \rf{ell(3)})
$$\ba{ll}
\ns\ds\BN^{\dbR^n}\big[\BN^D[h]\big](x)=\min_{\xi\in K}\(\BN^D[h](x+\xi)+\ell(T,x,\xi)\)\\
\ns\ds\qq=\min_{\xi\in K}\[\min_{\z\in K,x+\xi+\z\in\bar D}\(h(x+\xi+\z)+\ell(T,x+\xi,\z)+\ell(T,x,\xi)\)\]\\
\ns\ds\qq\ges\min_{\xi\in K}\[\min_{\z\in K,x+\xi+\z\in\bar D}\(h(x+\xi+\z)+\ell(T,x,\xi+\z)\)\]+\d_0=\BN^D[h](x)+\d_0.\ea$$
Also,
$$\BN^{\dbR^n}[h](x)=\min_{\xi\in K}\(h(x+\xi)+\ell(T,x,\xi)\)\les\min_{\xi\in K,x+\xi\in\bar D}\(h(x+\xi)+\ell(T,x,\xi)\)=\BN^D[h](x),\qq x\in\dbR^n.$$
We claim that under \rf{inf<}, the above equality holds. In fact, let
$$\min_{\xi\in K}\(h(x+\xi)+\ell(T,x,\xi)\)=h(x+\xi_0)+\ell(T,x,\xi_0),$$
for some $\xi_0\in K$. If $x+\xi_0\notin\bar D$, then by \rf{inf<}, the above leads to
$$\ba{ll}
\ns\ds\min_{\xi\in K}\(h(x+\xi)+\ell(T,x,\xi)\)=h(x+\xi_0)+\ell(T,x,\xi_0)\\
\ns\ds>h(x+\xi_0)>\min_{\xi\in K,x+\xi\in\bar D}\(h(x+\xi)+\ell(T,x,\xi)\),\ea$$
which is a contradiction. Hence, our claim holds and for $x\in\dbR^n\setminus\bar D$,
\bel{min=min}\ba{ll}
\ns\ds\min\Big\{h^D(x),\BN^{\dbR^n}[h^D](x)\Big\}=\min\Big\{\BN^D[h](x),
\BN^{\dbR^n}\big[\BN^D[h]\big](x)\Big\}=\BN^D[h](x);\ea\ee
For $x\in\bar D$,
$$\ba{ll}
\ns\ds\min\Big\{\min\big\{h(x),\BN^D[h](x)\big\},\BN^{\dbR^n}\[\min\big\{h,
\BN^D[h]\big\}\](x)\Big\}\\
\ns\ds=\min\Big\{h(x),\BN^D[h](x),\BN^{\dbR^n}[h](x),\BN^{\dbR^n}\big[
\BN^D[h]\big](x)\Big\}\\
\ns\ds=\min\Big\{h(x),\BN^{\dbR^n}[h](x)\Big\}=\min\Big\{h(x),\BN^D[h](x)\Big\}.\ea$$
These imply
\bel{}\big(h^D\big)^{\dbR^n}(x)=h^D(x),\qq x\in\dbR^n.\ee
From the above discussion, by a standard argument, we see that the value function $\wt V(\cd\,,\cd)$ of Problem ($\wt{\rm C}$), under some mild conditions, is the unique viscosity solution to the HJB QVI \rf{HJB-1}. We obtained
$$V(T,x)=\wt V(T,x),\qq\forall x\in\dbR^n.$$
But this does not mean that
\bel{V=V}V(t,x)=\wt V(t,x),\qq\forall(t,x)\in[0,T]\times\dbR^n.\ee
The reason is that it is unknown whether HJB QVI \rf{HJB-1} admits a unique viscosity solution in the function class that the value function $V(\cd\,,\cd)$ belongs to. On the other hand, if this were the case, then \rf{V=V} would hold. Consequently, under proper conditions that make $\wt V(\cd\,,\cd)$ to be globally Lipschitz (or H\"older) continuous, one should have the same continuity for the value function $V(\cd\,,\cd)$. However, from our discussion in Section 3, this seems to be unlikely. Hence, we end up with a challenging open questions:

\ms

\sl {\rm(i)} Is the viscosity solution to HJB QVI \rf{HJB-1} unique within the class of locally Lipschitz (or H\"older) continuous functions?

\ms

{\rm(ii)} Should there be some additional conditions for the HJB QVI so that it could characterize the value function?
\rm

\section{Concluding Remarks}

In this paper, we have introduced an intrinsic condition under which, together with other routine conditions, the value function of the optimal impulse control with terminal state constraint is continuous. This makes a big step towards the characterization of the value function. Due to the presence of the terminal state constraint, the value function is only locally Lipschitz (or H\"older) continuous and its growth is not slower than the impulse cost. Therefore, the available techniques are not enough to characterizing the value function as the unique viscosity solution of the HJB QVI. Moreover, efforts are made to the discussion on a seemingly equivalent optimal impulse control problem without terminal state constraint. It leads to a challenging open question about the uniqueness of viscosity solutions to the HJB QVI.

\baselineskip 18pt
\renewcommand{\baselinestretch}{1.2}

\end{document}